\documentclass[12pt]{article}
\usepackage{epsfig}
\title{Pentagram Spirals}
\author{Richard Evan Schwartz \footnote{Supported by N.S.F. grant DMS-0604426}}

\newtheorem{theorem}{Theorem}[section]

\newtheorem{lemma}[theorem]{Lemma}

\newtheorem{corollary}[theorem]{Corollary}
\newtheorem{conjecture}[theorem]{Conjecture}

\def\startproof{{\bf {\medskip}{\noindent}Proof: }}

\def\endproof{$\spadesuit$  \newline}

\def\P{\mbox{\boldmath{$P$}}}% 
\def\R{\mbox{\boldmath{$R$}}}% 
\def\T{\mbox{\boldmath{$T$}}}% 
\def\Z{\mbox{\boldmath{$Z$}}}% 

\begin{document}
\maketitle

\section{Introduction}

The {\it pentagram map\/} is a projectively
natural map defined on the space of
$n$-gons.  The case $n=5$ is classical;
it goes back at least to Clebsch
in the $19$th century and perhaps even
to Gauss.  Motzkin [{\bf Mot\/}] also
considered this case in $1945$.
I introduced the general version of the
pentagram map in $1991$.  See [{\bf Sch1\/}].
I subsequently published two additional
papers, [{\bf Sch1\/}] and [{\bf Sch2\/}],
on the topic.  Now there is a growing literature.
See the discussion below.

To define the pentagram map, one starts with
a polygon $P$ and produces a new polygon
$T(P)$, as shown at left in
Figure 1.1 for a convex hexagon.
As indicated at right, the map $P \to T^2(P)$
acts naturally on labeled polygons.  

\begin{center}
\resizebox{!}{2.2in}{\includegraphics{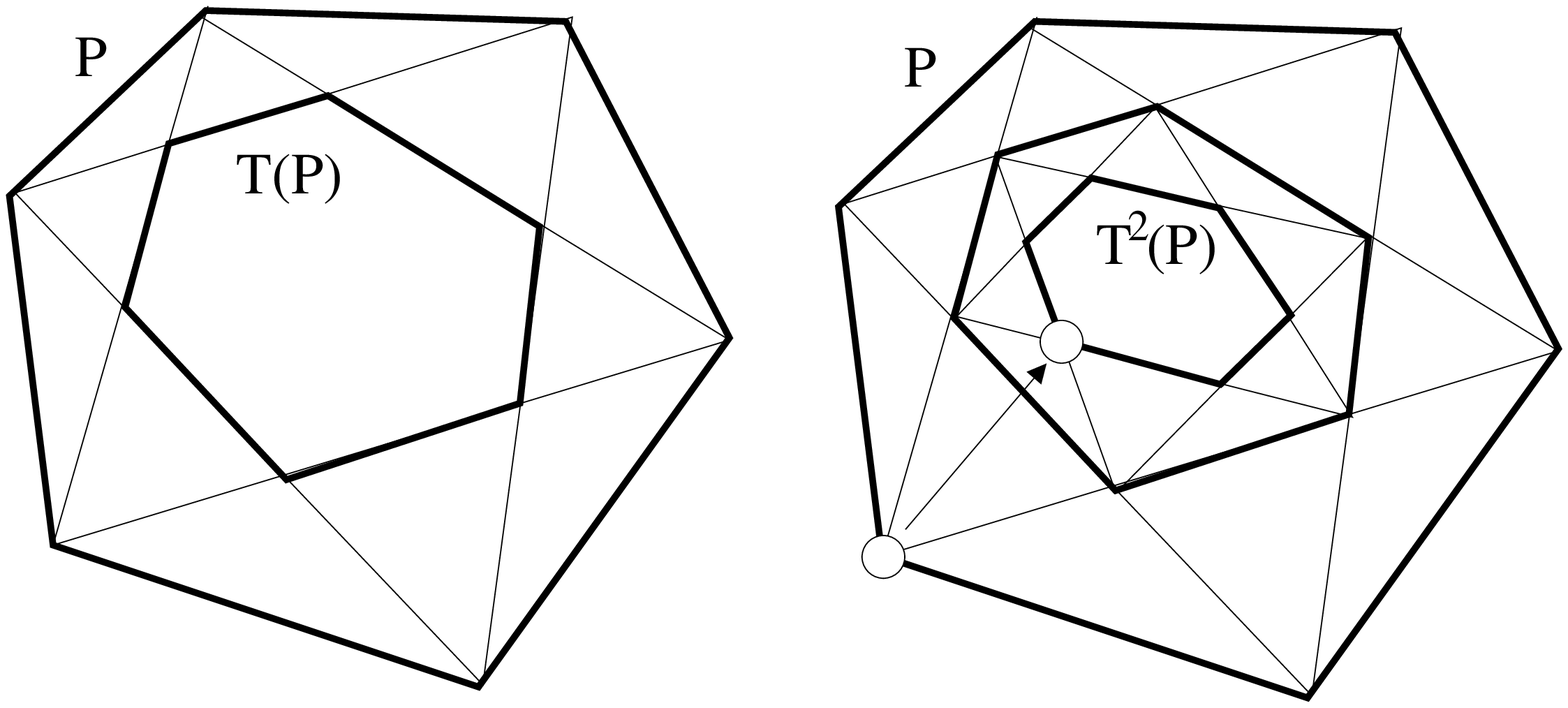}}
\newline
{\bf Figure 1.1:\/} The pentagram map
\end{center}

The pentagram
map is defined on polygons over any field.  More generally,
as I will discuss below, the pentagram map
is defined on the so-called twisted polygons.
The pentagram map commutes with projective
transformations and thereby induces a map
on spaces of projective equivalence classes
of polygons, both ordinary and twisted.

The purpose of this paper is to introduce a
variant of the pentagram map, which I will
call {\it pentagram spirals\/}.  The pentagram
spirals relate to the pentagram map much
in the way that logarithmic spirals
relate to circles.  I had the
idea for pentagram spirals many years ago,
but since there was not much interest in
the pentagram map, I decided not to pursue 
the idea.

\begin{center}
\resizebox{!}{5.2in}{\includegraphics{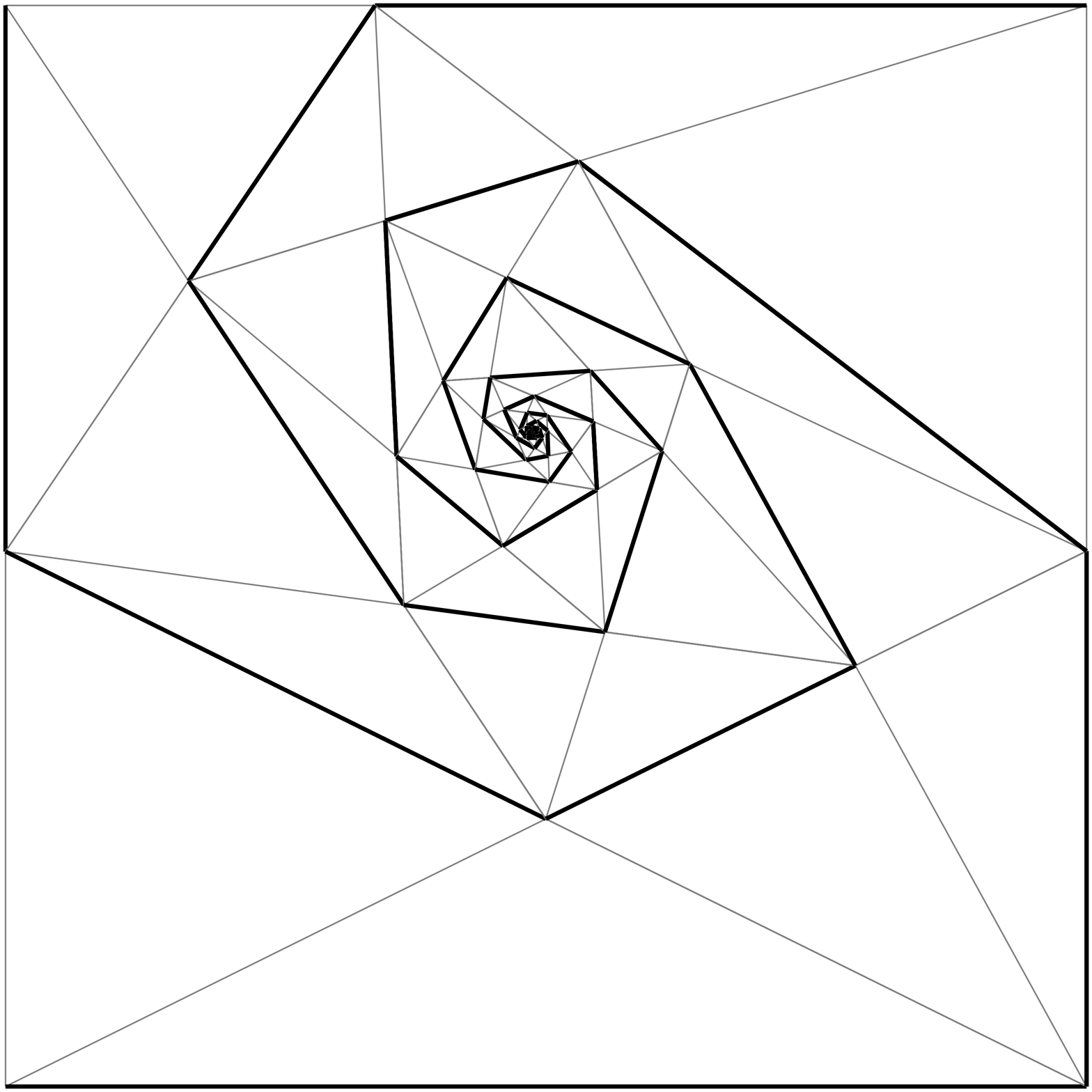}}
\newline
{\bf Figure 1.2:\/} The inward half of a pentagram spiral of type $(4,3)$.
\end{center}

In recent years, the pentagram map has
attracted a lot of attention,
thanks to the following developments.

\begin{enumerate}
\item In [{\bf Sch3\/}], I found a hierarchy of
integrals to the pentagram map, similar to the
KdV hierarchy.  I also related the pentagram map
to the octahedral recurrence, and observed
that the continuous limit of the pentagram map
is the classical Boussinesq equation.

\item In [{\bf OST1\/}], Ovsienko, Tabachnikov and I
showed that the pentagram map is a
completely integrable system when defined on
the space of projective classes of twisted polygons. We
also elaborated on the connection to the
Boussinesq equation.

\item In [{\bf Sol\/}] Soloviev 
showed that the pentagram map is completely
integrable, in the algebro-geometric sense,
on spaces of projective 
classes of real polygons and on spaces of
projective classes of complex polygons.
In particular Soloviev showed that
the pentagram map has a Lax pair and he
deduced the invariant Poisson structure from
the Phong-Krichever universal formula.

\item In [{\bf OST2\/}] (independently,
at roughly the same
time as [{\bf Sol\/}])
Ovsienko, Tabachnikov and I
showed that the pentagram map is a discrete,
completely integrable system, in the sense
of Liouville-Arnold, when defined on
the space of  projective classes of 
closed convex polygons.

\item In [{\bf Gli1\/}], Glick identified the pentagram
map with a specific cluster algebra, and found 
algebraic formulas for iterates of the map which are
similar in spirit to those found by Robbins and Rumsey
for the octahedral recurrence.

\item In [{\bf GSTV\/}], Gekhtman, Shapiro, Tabachnikov,
Vainshtein generalized the pentagram map to similar
maps using longer diagonals, and defined on spaces
of so-called {\it corrugated polygons\/} in higher
dimensions.  The work in [{\bf GSTV\/}] generalizes
Glick's cluster algebra.

\item In [{\bf MB1\/}], Mari-Beffa defines higher
dimensional generalizations of the pentagram map
and relates their continuous limits to
various families of integrable PDEs. See also
[{\bf MB2\/}]. 

\item In the recent [{\bf KS\/}],
Khesin and Soloviev obtain definitive results about
higher dimensional analogues of the pentagram map, 
their integrability, and their connection to 
KdV-type equations.

\item In the preprint [{\bf FM\/}], Fock and Marshakov
relate the pentagram map to, among other things,
Poisson Lie groups.

\item The preprint [{\bf KDiF\/}] discusses many aspects
of the octahedral recurrence, drawing connections to the
work in [{\bf GSTV\/}].

\end{enumerate}

Though this is not directly related to the
pentagram map, it
seems also worth mentioning the recent paper [{\bf GK\/}]
of Goncharov and Kenyon, who study a family of
cluster integrable systems.  These systems are
closely related to the octahedral recurrence which,
in turn, is closely related to the pentagram map.
\newline

Informally, a {\it pentagram spiral\/} is a bi-infinite
polygonal path $P$ in the projective plane such that
some finite power of 
the pentagram map carries $P$ to itself when $P$
is considered as an unlabeled path. \S 3.1 has
a formal definition. The global
combinatorics of how this is done allows one to
describe the {\it type\/} of the spiral by a
pair of integers $(n,k)$.  For instance, $k$
is the smallest integer such that 
$T^k(P)=P$, as an unlabeled path. The combinatorics
of the situation will be discussed in \S 3.2.

For every pair $(n,k)$ with $n \geq 4$ and $k=1,...,(n-1)$,
I will introduce a {\it pentagram spiral\/} of type
$(n,k)$. 
 I will focus on the case when the spirals
are what I call {\it properly locally convex\/},
or PLC for short.  The example in Figure 1.2 is PLC.
Though it is not nearly as obvious as in the case
of polygons, the basic constructions which generate
the polygon spirals just depend on drawing and
intersecting lines in the projective plane.  Thus,
they make sense over essentially any field.
However, we shall be interested mainly in the
PLC case.

We label the vertices of the spiral by consecutive
integers, so that the integers increase as the
spiral moves inwards.  A labeled pentagram spiral is
really just the same thing as a pentagram spiral with
a distingished vertex.  In \S 3 we will prove

\begin{theorem}
\label{cell}
The space ${\cal C\/}(n,k)$ of projective equivalence
classes of labeled PLC pentagram 
spirals of type $(n,k)$ has dimension 
$(2n-8)+k$ and is diffeomorphic to an open ball.
\end{theorem}

This result really amounts to describing how one
generates pictures like Figure 1.2.  The space
${\cal C\/}(n,k)$ should be seen as a relative
of the space ${\cal C\/}(n)$ of projective
classes of closed convex $n$-gons.  The space
${\cal C\/}(n)$ has dimension $2n-8$ and is
diffeomorphic to an open ball.

The spaces ${\cal C\/}(n,k)$ are naturally {\it shift spaces\/},
in that there is a natural map 
\begin{equation}
T_{n,k}: {\cal C\/}(n,k) \to {\cal C\/}(n,k)
\end{equation}
which just amounts to
moving the distinguished vertex inwards.
During the course of our proof of
Theorem \ref{cell}, we will define
$T_{n,k}$ from several points of view.

When properly interpreted, the map $T_{n,k}$ is a $d$th root
of the pentagram map, where $d=2k/(2n+k)$.
At the same time, when $n$ is large and $k$ is small, the
space ${\cal C\/}(n,k)$ is an approximation of the space
${\cal C\/}(n)$.  Thus,
$T_{n,k}$ in these cases is a very high root of
a map which is close to the pentagram map.  Independent
of the intrinsic beauty of the pentagram spirals, it seems
useful to have these high roots of approximations to the
pentagram map.

In \S 4 we will introduce projectively natural
coordinates on the space ${\cal C\/}(n,k)$ and
exhibit a $T_{n,k}$ invarant function.  This
is the analog of the invariant function
used in [{\bf Sch1\/}] and [{\bf Sch2\/}].  For
experts, our first invariant function is the
analogue of what is, in the closed case,
one of the Casimirs for the invariant Poisson
structure.

In \S 5 we will prove the following compactness result,
which is similar in spirit to a similar result in
[{\bf Sch1\/}].
\begin{theorem}
\label{compact}
The orbits of $T_{n,k}$ have compact closure in ${\cal C\/}(n,k)$.
\end{theorem}
Theorem \ref{compact} says that up to projective
equivalence one sees roughly the same
shape, over and over again, as one moves inwards
or outwards along the spirals.  The proof just
involves showing that our invariant function has
compact level sets.

In \S 6 we will use Theorem \ref{compact}
to deduce several geometric corollaries
about PLC pentagram spirals.

\begin{theorem}
\label{omega}
The $\omega$-limit set of a PLC pentagram spiral,
in the projective plane, is a union
of a single point and a single line.  Equivalently, the
support of the triangulation associated to a pentagram spiral is
projectively equivalent to a punctured Euclidean plane.
\end{theorem}

Theorem \ref{omega} implies that the forward direction
of $P$ spirals down to a single {\it limit point\/}.
One might wonder about the nature of this spiraling.

\begin{theorem}
\label{wind}
A PLC pentagram spiral winds infinitely many
times around its limit point.
\end{theorem}

Theorems \ref{omega} and \ref{wind} pin down
some of the rough geometry 
of the triangulations associated to the
pentagram spirals. Theorem \ref{omega}
is the analog
of the result in [{\bf Sch1\/}] which says
that the pentagram map shrinks arbitrary
convex polygons to single points.
The triangulations associated to the pentagram
spirals are locally
the same as the triangulations one sees
when one takes the full orbit of a convex polygon
under the pentagram map. However, the global structure
of the tilings is different.  

The main purpose of this paper is geometric.
The above results answer probably the most basic
geometric questions one would want to know about
PLC pentagram spirals.  I believe that the
deep algebraic structure underlying the pentagram
map is also present in the pentagram spirals, and
I hope that this paper imspires future work on
these objects.
In the informal \S 7, I will discuss some 
computer experiments, conjectures, and
topics for further study.

I wrote a Java program which allows the user
to draw the pentagram spirals for smallish
values of $n$ and $k$, and also to watch the
spirals evolve under the map $T_{n,k}$.  
You can download this program at 
\newline
\newline
{\bf http://www.math.brown.edu/$\sim$res/Java/SPIRAL.tar\/}
\newline
\newline
I strongly suggest that the reader interested in this
paper download the program and play with it.  I think
that the program greatly enriches the paper.

This paper contains rigorous proofs
of all the main results, but it seems
worth mentioning that I checked the
results numerically using my computer
program.  For instance, my program
draws the pentagram spirals using
exactly the formula for $T_{n,k}$ given
in \S \ref{seedmap1}.  The program
also lets the user see that the invariant
function $Z$ defined in \S \ref{first}
is indeed an invariant of the map
$T_{n,k}$.

I first thought about the pentagram spirals many
years ago when I was
talking to Peter Doyle about his so-called Doyle spirals.
The Doyle spirals are
 circle packings which relate to the hexagonal circle
packing much in the way that the pentagram spirals
relate to the pentagram map.  Recently I walked 
through the streets of Nice on a pleasant spring night
and thought back to the puzzles of my youth. In a
nostalgic mood, I decided to send a short note about the pentagram
spirals to some of the researchers interested in the
pentagram map.  After sending the note, I got the bug
again and decided
to write a more systematic paper about them.

I'd like to thank Valentin Ovsienko and Sergei Tabachnikov,
my usual collaborators on the pentagram map, for many
discussions about the pentagram map and related areas
of mathematics.   This work was carried 
out during my sabbatical at Oxford in
2012-13.  I would especially like to
thank All Souls College, Oxford, for
providing a wonderful research environment.
My sabbatical was funded from many sources.
I would like to thank the National
Science Foundation, All Souls College,
the Oxford Maths Institute,
the Simons Foundation, the Leverhulme Trust, the
Chancellor's Professorship, and
Brown University for their support during this time
period.

\newpage

\section{Preliminaries}

\subsection{Projective Geometry}

The {\it real projective plane\/}
$\R\P^2$ is the space of lines
through the origin in $\R^3$.
Such lines are denoted by
$[x:y:z]$.  This is the line
consisting of all vectors
of the form $(rx,ry,rz)$ with
$r \in \R$.

In the usual way, we think of
$\R^2$ as an affine patch of
$\R\P^2$.  Concretely,
the inclusion is given by
\begin{equation}
(x,y) \to [x:y:1].
\end{equation}

A {\it line\/} in $\R\P^2$ is a collection
of points represented by all the lines in
a plane through the origin in $\R^3$.
The set $\R\P^2-\R^2$ is a single line,
called the {\it line at infinity\/}.
All other lines in $\R\P^2$ intersect
$\R^2$ in a straight line.  Conversely,
any straight line construction
in $\R^2$ extends naturally to a
straight line construction in
$\R\P^2$.  When we make our constructions,
we will draw things in the plane (of course)
but we really mean to make the constructions
in the projective plane.

A {\it projective transformation\/}
is a self-homeomorphism of $\R\P^2$ induced
by the action of an invertible linear
transformation.  Projective transformations
permute the lines of $\R\P^2$ and are
in fact analytic diffeomorphisms.
Conversely any homeomorphism of
$\R\P^2$ which carries lines to lines
is a projective transformation.

A {\it projective construction\/} is one in which
points and lines are produced by the following
two operations:
\begin{itemize}
\item Given distinct points $a$ and $b$, take
the line $(ab)$ through $a$ and $b$.
\item Given distinct lines $l$ and $m$, take
the intersection $l \cap m$.
\end{itemize}
We use the notation $(ab)(cd)$ to denote the
intersection of the line $(ab)$ with the
line $(cd)$.  Since projective transformations
carry lines to lines, any projective
construction commutes with the action
of the group of projective transformations.

A subset of $\R\P^2$ is {\it convex\/} if
some image of that subset under a projective
transformation is a convex subset of
$\R^2$ in the ordinary sense.  For instance,
a hyperbola in the plane extends to a closed
loop in $\R\P^2$ which bounds a convex subset
on one side.  Proof: one can move the hyperbola by
a projective transformation so that it is a circle
in the plane.

\subsection{The Cross Ratio}

The {\it inverse cross ratio\/} of $4$ real
numbers $a,b,c,d \in \R$ is the quantity
\begin{equation}
[a,b,c,d]=\frac{(a-b)(c-d)}{(a-c)(b-d)}
\end{equation}
When $a<b<c<d$, the quantity
$[a,b,c,d]$ lies in $(0,1)$.
We will usually consider this situation.

Given $4$ collinear points
$A,B,C,D$ in the projective plane, we choose
some projective transformation which identifies
these points with $4$ numbers on the
$x$-axis, and then we take the cross ratio
of the first coordinates of these numbers.
This lets us define $[A,B,C,D]$.  The
result is independent of any choices made.
In particular
\begin{equation}
[A,B,C,D]=[T(A),T(B),T(C),T(D)]
\end{equation}
for any projective transformation $T$.

\subsection{The Corner Invariants}
\label{invariant}

In [{\bf Sch3\/}] I introduced the notion
of {\it corner invariants\/} of a polygon
in the projective plane.  We also used these
invariants in
[{\bf OST1\/}] and [{\bf OST2\/}].
For the basic definition we will follow
the notation in [{\bf OST2\/}], but then
we will revert back to an interpretation
of the invariants given in [{\bf Sch3\/}].

The {\it corner invariants\/} of a polygonal
path $P$ with successive vertices $\{v_i\}$
are defined as follows.

We define
\begin{equation}
x_{3}=
[v_0,v_1,(v_0v_1)(v_2v_3),(v_0v_{1})(v_3v_{4})]
\end{equation}

\begin{equation}
x_4=[v_4,v_3,(v_4v_3)(v_2v_1),(v_4v_3)(v_1v_0)].
\end{equation}

The remaining invariants are defined by shifting
the indices by $2k$ for $k \in \Z$.  Figure 2.1
shows a picture of the construction.

\begin{center}
\resizebox{!}{3in}{\includegraphics{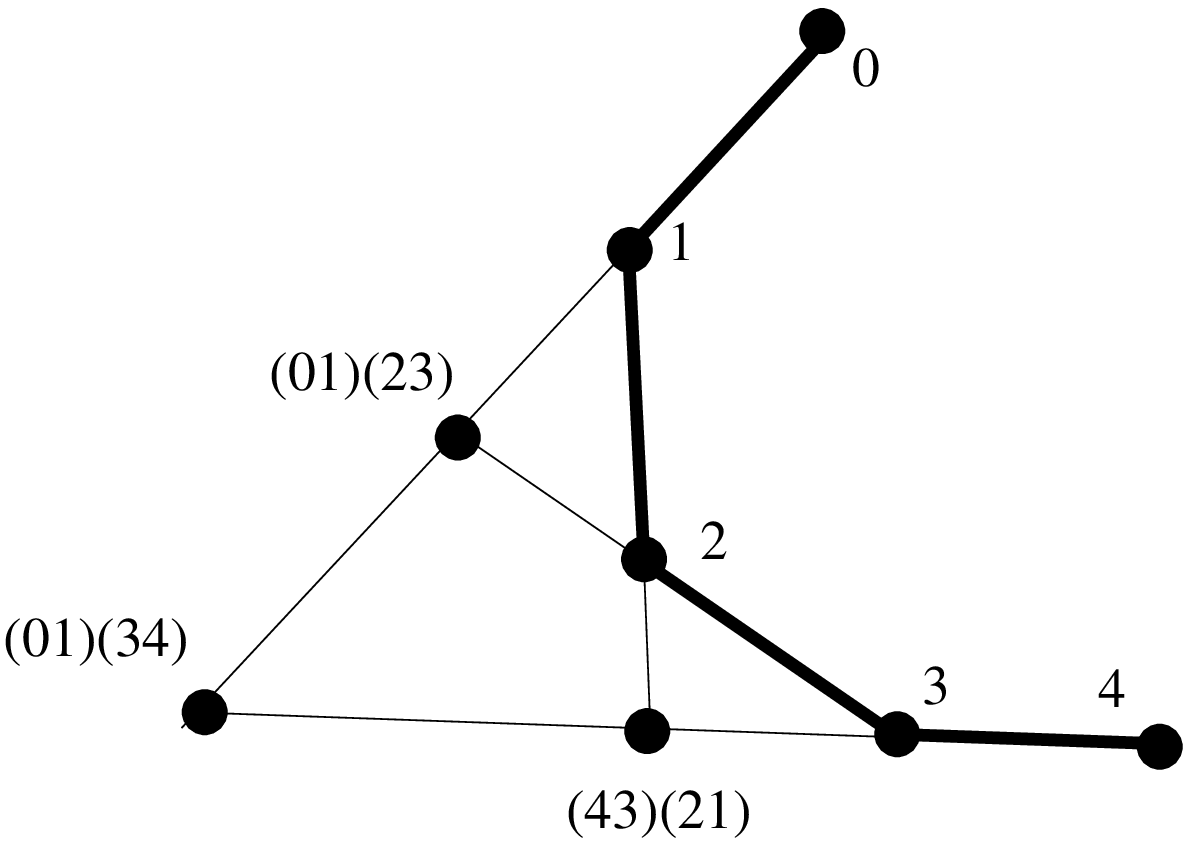}}
\newline
{\bf Figure 2.1:\/} The corner invariants
\end{center}

When the curve is such that every $5$ consecutive points
make the vertices of a convex pentagon, the invariants
all lie in $(0,1)$.

The corner invariants are projectively natural.
They provide projectively natural coordinates
on the space of polygons.  Two polygons are projective
equivalent if and only if they have the same corner
invariants.  See [{\bf Sch3\/}] for details.
We can express the pentagram map in these coordinates,
provided that we ``break symmetry'' and choose a less
than canonical labeling scheme for $P$ and $T(P)$.
We label the vertices of $P'=T(P)$ so that
\begin{equation}
v_1'=(v_0v_2)(v_1v_3).
\end{equation}
The remaining labels are obtained by shifting the
indices.  In [{\bf OST2\/}], we called this the
{\it right convention\/}.
Using the right convention, we have

\begin{equation}
\label{penta}
x'_2=x_4 \frac{1-x_5x_6}{1-x_1x_2}, \hskip 30 pt
x'_3=x_3 \frac{1-x_1x_2}{1-x_5x_6}.
\end{equation}
The remaining equations are obtained by shifting
the indices by $2k$.  
These equations are less than perfectly
symmetric, on account of the symmetry-breaking
labeling convention.  However, they served our
purposes in [{\bf OST1\/}] and [{\bf OST2\/}].
In the next section we will explain a more 
symmetric picture.

\subsection{Tiling Interpretation of the Coordinates}
\label{tiling}

For convenience we work with bi-infinite paths.
Before getting to the paths, however, we will
dicsuss a seemingly different construction.
We consider a hexagonal tiling of the plane,
by right-angled isosceles triangles whose hypotenuse
is horizontal.  Figure 2.2 shows a small part
of this triangulation.

\begin{center}
\resizebox{!}{3in}{\includegraphics{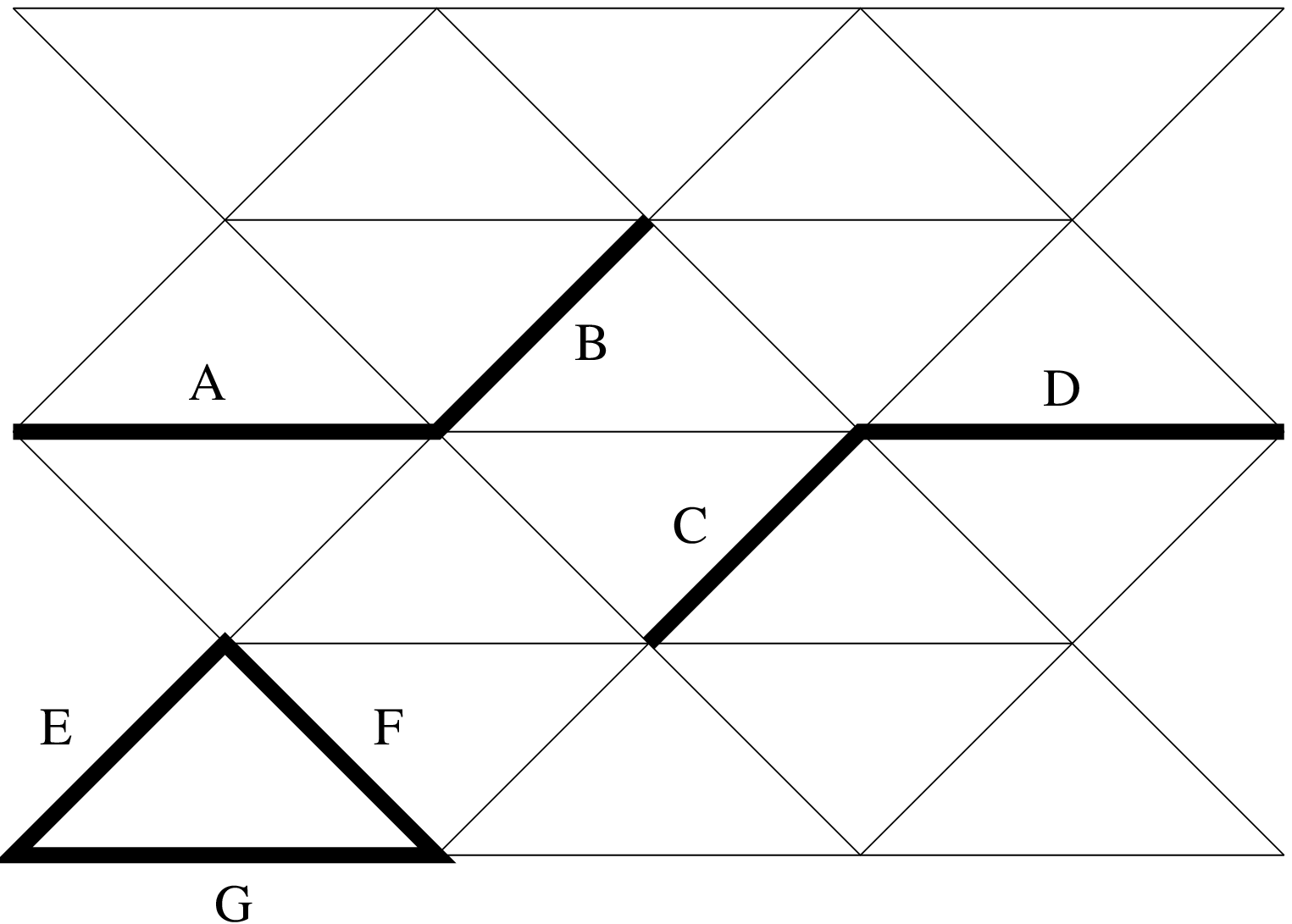}}
\newline
{\bf Figure 2.2:\/} The corner invariants
\end{center}

We label the edges of the triangulation by real
numbers and we insist on the compatibility relations
\begin{equation}
AB=CD, \hskip 30 pt G=1-EF.
\end{equation}
These relations are meant to hold for the labels
of all the isometric images of these configurations.
The configurations may be translated or reflected.
We call such a labeling a {\it pentagram labeling\/}.

If we label the diagonal edges of one row of
diagonals, we can fill in the remaining labels using
the compatibility rules.  This will certainly remind
some readers of the octahedral recurrence, and indeed
we worked out the connection explicitly in 
[{\bf Sch3\/}].

To see the connection between the pentagram labelings
and the pentagram map, we fill in two rows of diagonal
edges with the corner invariants, as in Figure 2.3.

\begin{center}
\resizebox{!}{1.7in}{\includegraphics{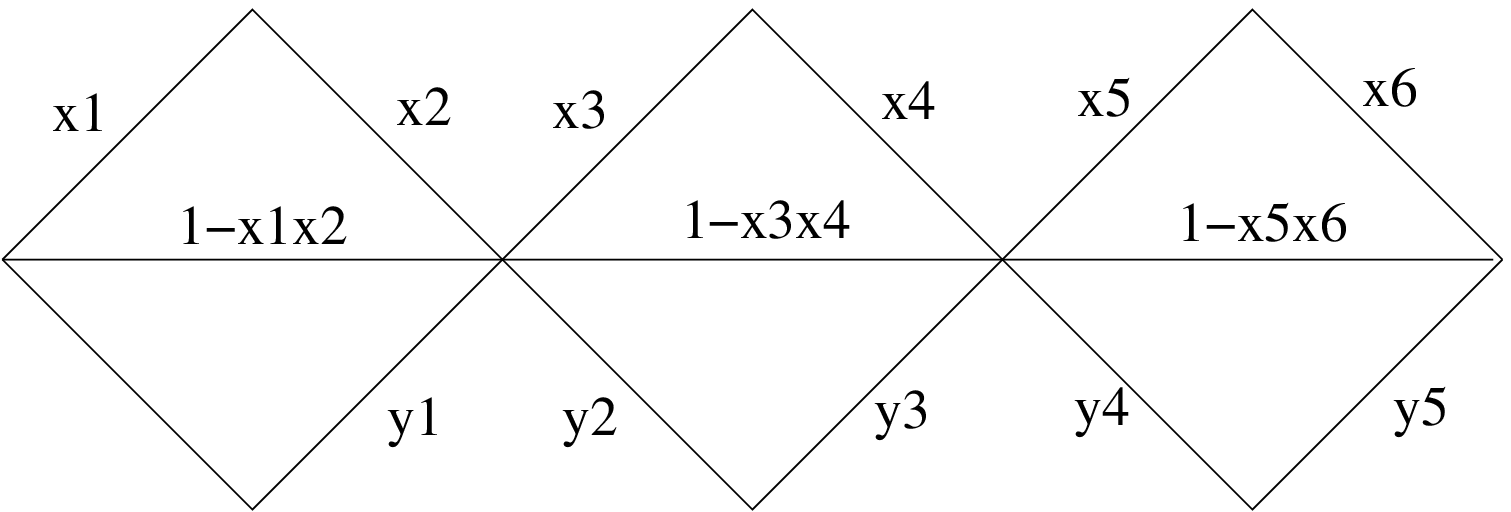}}
\newline
{\bf Figure 2.3:\/} The pentagram map seen on the tiling
\end{center}

For ease of labeling we set $y_i=x'_i$.  The
compatibility rules above express the
$y$-coordinates in terms of the $x$-coordinates,
and the formula is exactly as given for the pentagram
map in Equation \ref{penta}.
The cost of breaking symmetry 
is that the correspondence between the triangulation
labels with the corner invariants is somewhat
asymmetric.  Indeed, if we wanted to continue this
correspondence down to the next row of the tiling, so
to speak, we would have to switch from the right to
the left labeling convention.

In [{\bf Sch3\/}] we used
a different scheme, whereby the corner invariants
corresponded not to vertices of the polygon but rather
flags. The system in [{\bf Sch3\/}] worked perfectly
from the point of view of lining up the pentagram 
picture with the tiling picture, but
the apparatus was somewhat cumbersome.
Ultimately, we dropped this scheme in
[{\bf OST1\/}] and [{\bf OST2\/}],
settling on something less canonical but more
businesslike.  In \S 4 we will try for the
best of both worlds, choosing a convention
for the corner invariants which captures the
symmetry in [{\bf Sch3\/}] but retains the
efficient nature of [{\bf OST1\/}] and [{\bf OST2\/}]

In case we are working in ${\cal C\/}(n)$, the
labelings we get are periodic with respect to
a horizontal translation by $n$.  (We think of the
horizontal edges as having unit length.)  In this
case, the reader can probably see that the
products
\begin{equation}
E(P)=\prod_{i=1}^n x_{2i}, \hskip 30 pt
O(P)=\prod_{i=1}^n x_{2i+1}
\end{equation}
are such that
\begin{equation}
E(T(P))=O(P), \hskip 30 pt
O(T(P))=E(P).
\end{equation}
Thus $E$ and $O$ are invariants of the
square of the pentagram map.  
These are the first invariants of many.

\newpage

\section{The Space of Pentagram Spirals}
\subsection{Basic Definitions}
\label{maindef}

It is actually not so easy to give a formal
definition of a pentagram spiral.  We will
start out with an easy provisional definition
but then we will explain the problems with it.
Our final definition is somewhat more technical.
In the discussion, $T$ denotes the pentagram map,
defined on unlabeled paths.
\newline
\newline
{\bf Definition:\/}
A bi-infinite polygonal path $P \subset \R\P^2$ is a
{\it weak pentagram spiral\/} if
path $P \subset \R\P^2$ with the following
there is some integer $k>0$ such that
$T^k$ is defined on $P$ and
$T^k(P)=P$.  The smallest $k$ with this property
is called the {\it order\/} of $P$.
\newline
\newline
It is important to note that perhaps $T^k(P)=P$ only
in the unlabeled sense.  For instance, when $k$ is
even, there is a canonical notion of an action of
$T$ on labeled bi-infinite paths.  In this
situation, it might be the case that $T^k(P)$
and $P$ do not agree as labeled paths.

The definition above is too broad to be of use to
us. For example, suppose we have a
polygon $Q$ which is periodic under the pentagram map,
in the literal sense that $T^k(Q)=Q$. The pentagram
map has many periodic points when acting on projective
classes of polygons, but here we mean that the
actual polygon is periodic with respect to
the pentagram map. We do not
have an explicit example of this, but presumably
it can happen in the non-convex case.  Then
we could take $P$ to be a bi-infinite path which
winds around $Q$ infinitely often in both
directions. The path $Q$ would be a weak pentagram
spiral according to the definition above.

Here is the definition we care about.
\newline
\newline
{\bf Definition:\/} A {\it propertly locally convex pentagram spiral\/}
is a weak pentagram spiral $P$ of
order $k$ such that the iterates $T^j(P)$ are embedded,
locally convex, and pairwise disjoint, for $j=0,...,(k-1)$.
\newline

We will usually abbreviate {\it properly locally convex\/}
to {\it PLC\/}.
Figure 3.1 shows an example of a PLC pentagram spiral
of order $2$.  These spirals are meant to go outwards
as well as inwards.

\begin{center}
\resizebox{!}{3.8in}{\includegraphics{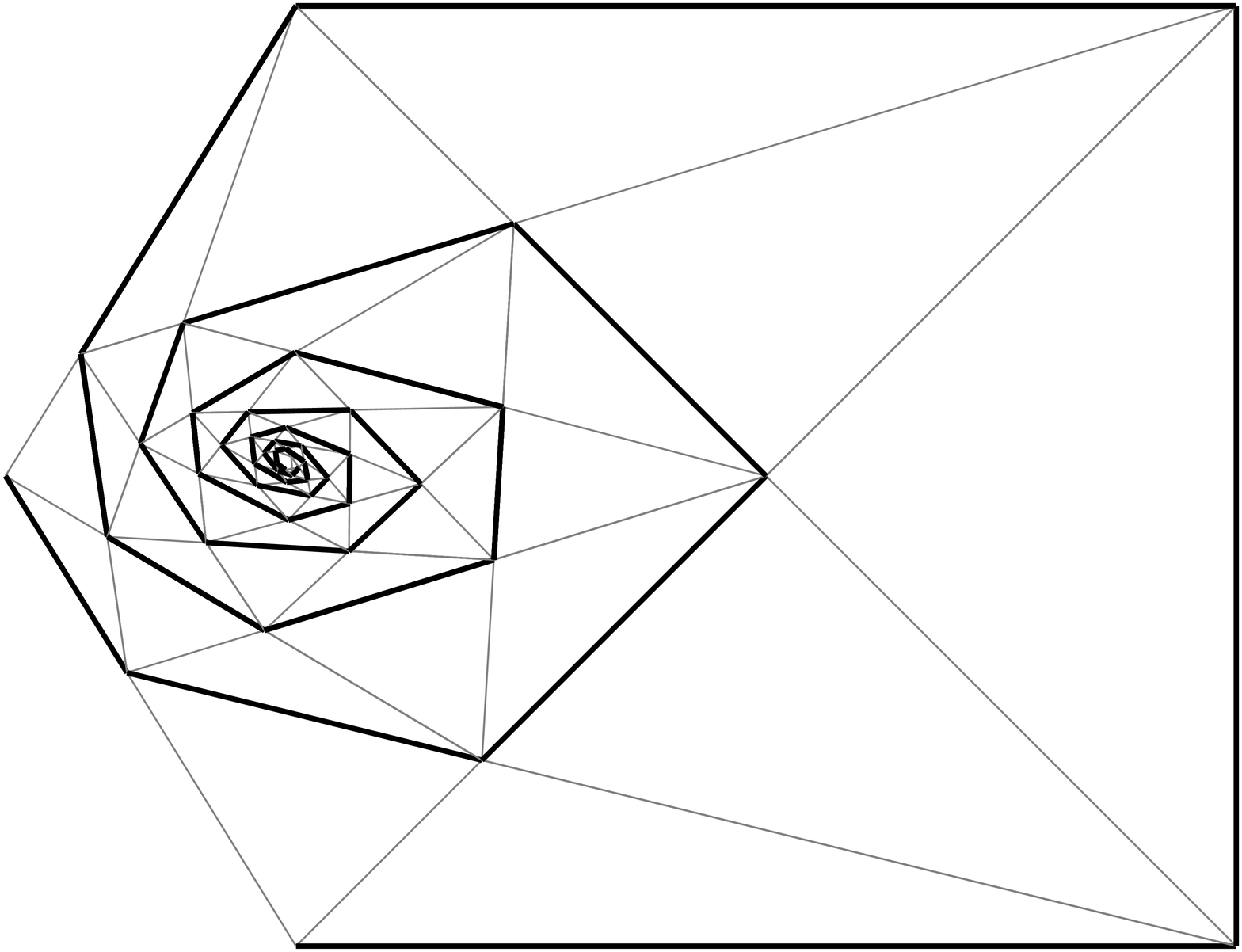}}
\newline
{\bf Figure 3.1:\/} A PLC pentagram spiral of type $(5,2)$.
\end{center}

We are not interested in pentagram spirals which are not
PLC, but we give a formal definition for the interested
reader.  To make this definition go more smoothly, we will
assume Theorem \ref{omega}.  Suitably normalized,
a PLC pentagram spiral defines a triangulation
of the punctured plane.  Figure 3.1 shows the
``inner half'' of this triangulation.
We call this tiling a PLC tiling.

Let $\tau$ denote a PLC tiling.  An {\it adapted immersion\/}
of $\tau$ is continuous map from $\tau$ into $\R\P^2$
which maps each line segment of $\tau$ to a line segment
in $\R\P^2$.  We also require the map to carry each
triangle of $\tau$ to a nontrivial triangle.
\newline
\newline
{\bf Definition:\/} A {\it pentagram spiral\/} is the
image of a PLC pentagram spiral under an adapted
immersion of the corresponding PLC tiling.
\newline

Essentially, a pentagram spiral is an object
which has the same locally combinatorial structure
as a PLC pentagram spiral. Since we only discuss
PLC pentagram spirals in this paper, the reader
need not absorb this last definition in order
to understand the rest of the paper.

\subsection{Combinatorics of Pentagram Spirals}

We will always work with PLC pentagram spirals,
though what we say usually works for a general
pentagram spiral.
We label the vertices of each pentagram spiral by
the integers, so that the numbers increase as one
moves inwards along the spiral.  In our example
in Figure 3.1, the map $T^2$ carries $P$ to itself
but $T^2(P_j)=P_{j+6}$ for all $j \in \Z$.
Here $P_j$ is the $j$th point of $P$.
In general, there is some number $\mu$, either
a whole or a half integer, such that
\begin{equation}
T^{2k}(P_j)=P_{j+2\mu}.
\end{equation}
The pair of numbers $(\mu,k)$ characterizes 
the combinatorial type of the tiling produced
by iterating $T$ on the spiral.  In 
Figure 3.1, we have $\mu=6$.  In Figure 1.2,
we have $\mu=5\frac{1}{2}$.

We find it more convenient to replace $\mu$ with
another invariant which captures the same
information.  We will use the pair
$(n,k)$ where $k$ is the order of the spiral
and
\begin{equation}
n=\mu-\frac{k}{2}.
\end{equation}
Geometrically, it turns out that $n$
counts the number of sides of the
``seed'' which generates the spiral.
We will explain precisely what we mean by a
seed in the next section, but informally,
the seed is the outer polygon in Figures 1.2 and
3.1. So, we have $n=4$ in Figure 1.2 and
$n=5$ in Figure 3.1.  

There is a natural map on the space of
pentagram spirals of type $(n,k)$, which we
call $T_{n,k}$.  The map $T_{n,k}$ shifts
the labeling of a spiral by one unit.  Thus
$T_{n,k}(P)$ is the same unlabeled spiral as
$P$, but the $k$th vertex of $T_{n,k}(P)$ is
the $(k+1)$st vertex of $P$.  It seems
at first that the map $T_{n,k}$ is trivial,
but in fact, on labeled pentagram spirals we have
\begin{equation}
T_{n,k}^{2n+k}(P)=T^{2k}(P).
\end{equation}

In the next several sections, we will develop the
idea of generating a PLC pentagram spiral from
a seed. The reason we do this is $2$-fold. 
First, the seeds give a convenient way
for drawing the spirals. 
Second, it turns out that $T_{n,k}$ can
be interpreted as a kind of evolution operator
on the set of seeds.  The spirals are
generated by considering the orbit of the
seed under powers of $T_{n,k}$.

Once we have the
basic constructions involving seeds, we will
prove Theorem \ref{cell}.  The basic idea
is to show that every seed generates a
PCL spiral and that every PLC spiral comes from
a seed. Finally, we will identify the version
of $T_{n,k}$ given in terms of seeds with the
shift map discussed above. 

\subsection{Seeds}

A {\it seed\/} of type $(n,k)$ is a 
strictly convex $n$-gon
with an additional point chosen in the
interior of each of
the last $k$ edges.  More precisely,
the vertices are
points $A_1,...,A_{n}$ and the additional
$B_{n-k+1},...,B_n$.  Here
$B_j$ lies in the interior of the
edge $A_jA_{j+1}$, with indices
taken mod $n$.  
We will sometimes denote our seeds as
$(A,B)$, where $A$ is short
for $\{A_1,...,A_n\}$ and
$B$ is short for $\{B_{n-k+1},...,B_n\}$.

\begin{center}
\resizebox{!}{3in}{\includegraphics{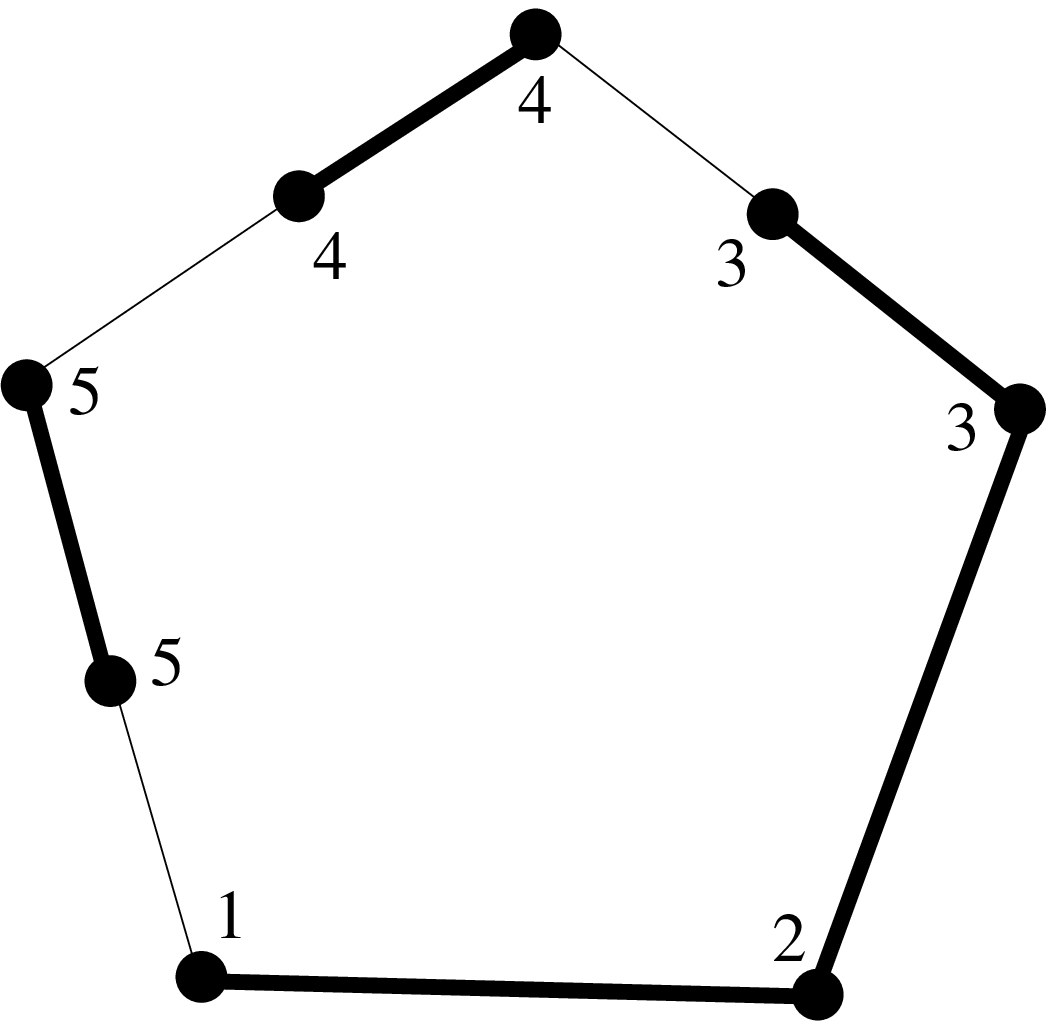}}
\newline
{\bf Figure 3.2:\/} A seed of type $(5,3)$.
\end{center}

Figure 3.2 shows an 
example of a seed of type $(5,3)$.
The labeling should be fairly obvious:
The vertex labeled $m$ denotes $A_m$
and the point labeled $m$ in the middle
of an edge denotes $B_m$.
We decorate the polygon with thick segments
of the following kind:
\begin{itemize}
\item $(A_iA_{i+1})$ for $i=1,...,(n-k+1)$.
\item $(A_i,B_i)$ for $i=(n-k+1),..,n$.
\end{itemize}
We call these segments the {\it spiral segments\/}.
The idea is that a PLC pentagram spiral and
its iterates under the pentagram map will
turn out to be an infinite union of spiral
segments, taken from an infinite union of
seeds which piece together in a way that we
explain in the next section.

\subsection{The Seed Map}
\label{seedmap1}

Figure 3.3 shows a special case of the general
projective construction which produces a
new seed from an old one.

\begin{center}
\resizebox{!}{3.4in}{\includegraphics{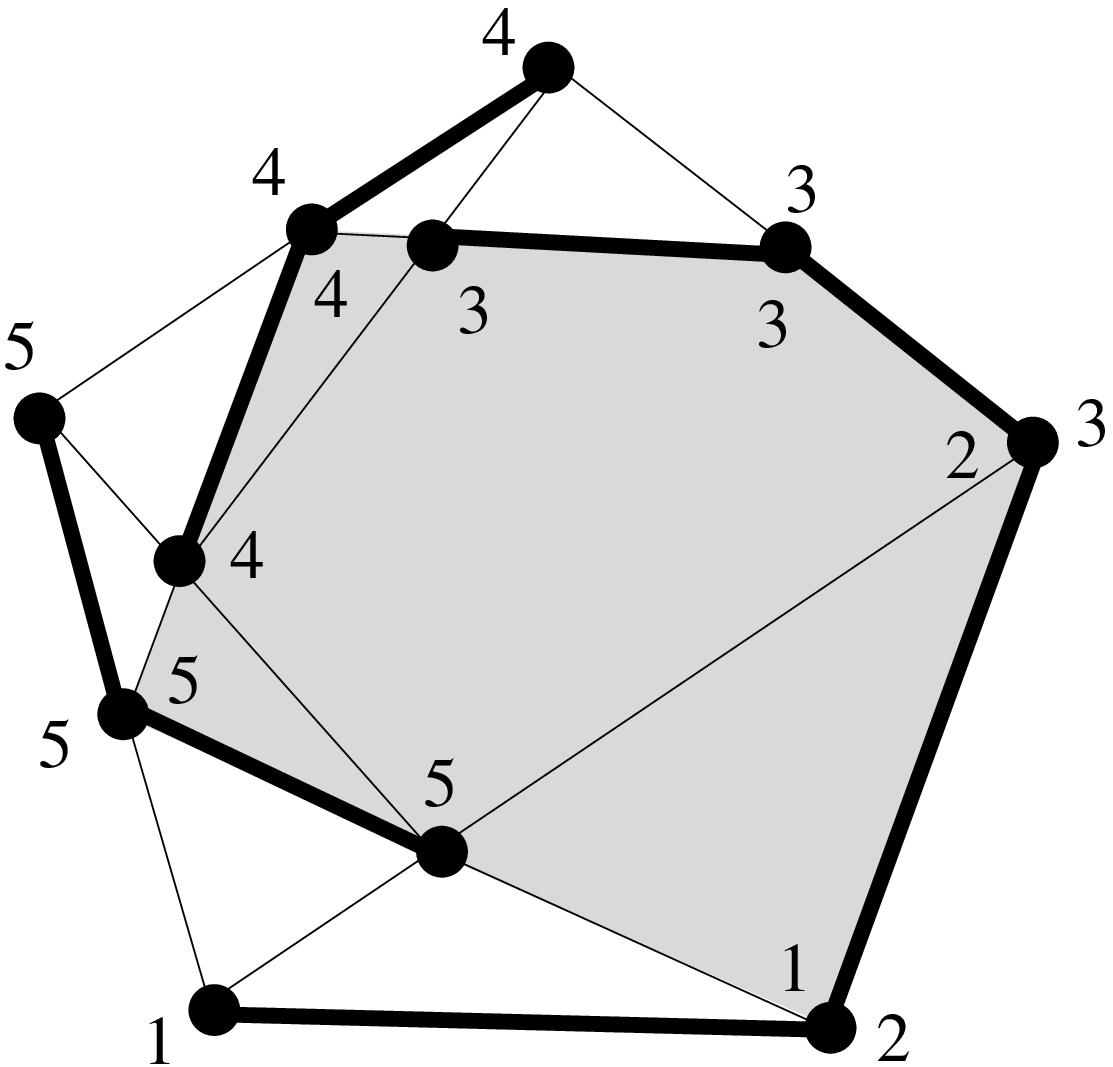}}
\newline
{\bf Figure 3.3:\/} One seed produces another.
\end{center}

The boundary
of the shaded region is the polygon supporting
the new seed.  The outside
numbers are the old seed labels and the inside
numbers are the new seed labels.
Putting a star for the new points, we have
\begin{itemize}
\item $A_i^*=A_{i+1}$ for $i=1,2$.
\item $A_i=B_i^*$ for $i=3,4,5$.
\item $B_5^*=(A_1A_2^*)(A_5^*A_1^*)$.
\item $B_i^*=(A_{i+1}B_{i+1}^*)(A_i^*A_{i+1}^*)$ for
$i=3,4$.
\end{itemize}
In the last item, we must
construct the point with the larger index
$(i=5)$ first.
Notice that the spiral segments on the two seeds
line up to produce what promises to be a union
of spirals.  Figure 3.6 below adds another seed, making
the connection to the pentagram map somewhat clearer.

\begin{center}
\resizebox{!}{2.5in}{\includegraphics{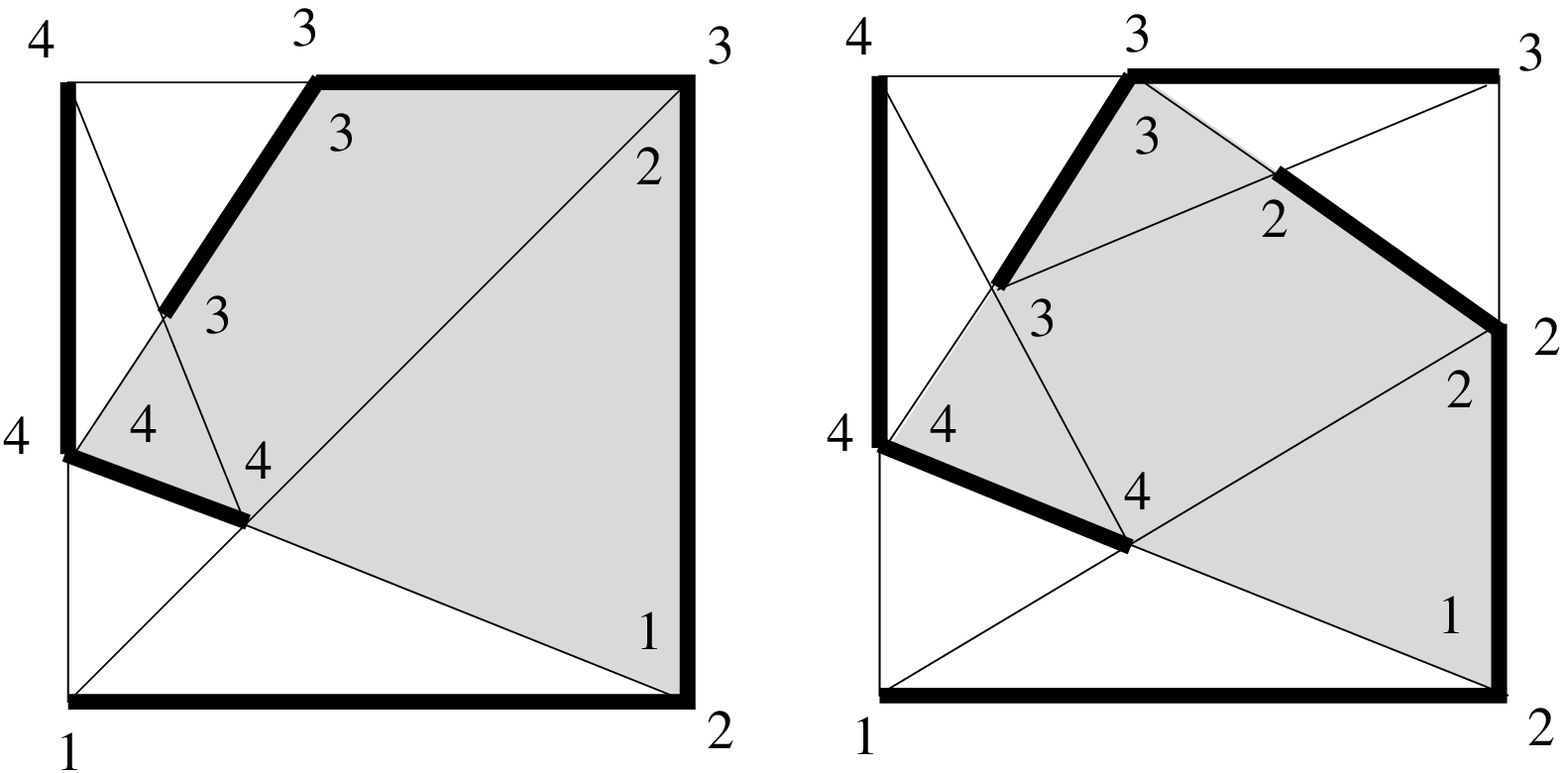}}
\newline
{\bf Figure 3.4:\/} 
The $(4,2)$ and $(4,3)$ cases of the construction.
\end{center}
Figure 3.4 shows two more examples.

In general, we define
\begin{equation}
T_{n,k}(A_1,...,A_n;B_{n-k+1},...,B_n)=
(A_1^*,...,A_n^*;B_{n-k+1}^*,...,B_n^*)
\end{equation}
according to the following rules.
\begin{itemize}
\item $A_i^*=A_{i+1}$ for $i=1,...,(n-k)$.
\item $A_i^*=B_i$ for $i=(n-k+1),...,n$.
\item $B_n^*=(A_1A_2^*)(A_n^*A_1^*)$.
\item $B_j^*=(A_{j+1}B_{j+1}^*)(A_j^*A_{j+1}^*)$ for
$j=(n-1),...,(n-k+1)$.
\end{itemize}
In the last item, it is important that the points are 
constructed going from the largest to the smallest
index, so that the map is a projective construction.

\begin{lemma} The starred points form a new
seed of the same type.
\end{lemma}

\startproof
The $A^*$ polygon is obtained from the
$A$ polygon by cutting off $k$ corners
using non-overlapping line segments.
Hence $A^*$ is a strictly convex polygon.
The triangle
$A_1A_1^*A_n^*$ is nondegenerate and oriented
counterclockwise. A routine
induction argument shows that
the triangles $A_jA_j^*B_j^*$ are also
nondegenerate and oriented counterclockwise for
$j=n,...,(n-k+2)$.  These facts, together with the
definition of the point, imply that each
$B_j^*$ lies in the interior of the segment 
$A_j^*A_{j+1}^*$
\endproof

\subsection{The Inverse Map}
\label{seedmap2}

Figure 3.5, which is a repeat of
Figure 3.3 but with different shading,
can be interpreted instead as an illustration of
how one can derive outer seed from the
inner one.  

\begin{center}
\resizebox{!}{3.5in}{\includegraphics{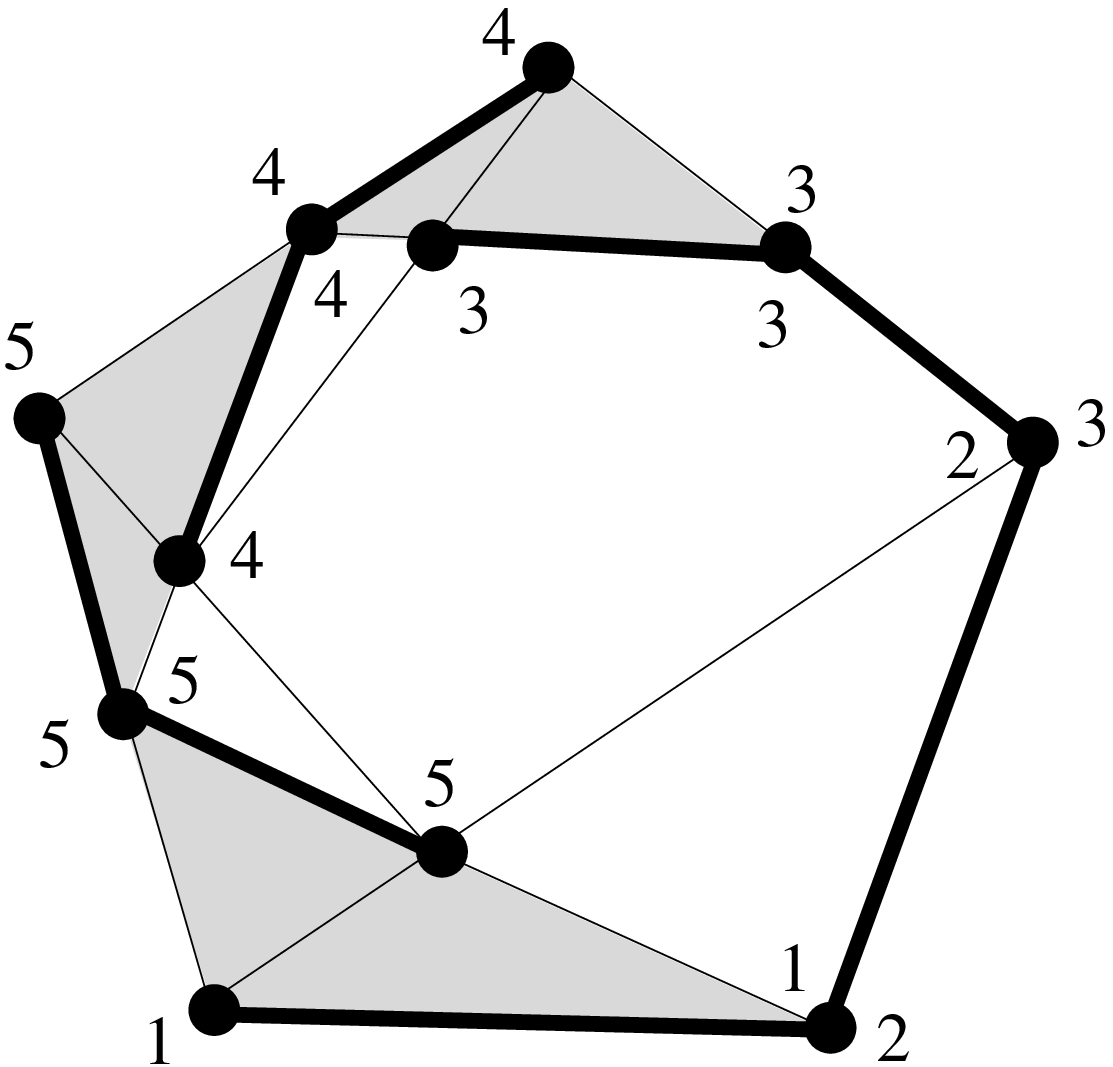}}
\newline
{\bf Figure 3.5:\/} Figure 3.3 repeated.
\end{center}

We leave it to the reader to check that, in general,
\begin{itemize}
\item $A_i=A_{i-1}^*$ for $i=2,...,(n-k+1)$.
\item $B_i=A_i^*$ for $i=(n-k+1),...,n$.
\item $A_i=(A_{i-1}B_{i-1})(B_{i-1}^*B_i^*)$ for $i=(n-k+2),...,n$.
\item $A_1=(A_nB_n)(A_2^*B_n^*)$.
\end{itemize}
This gives us a formula for the inverse map
$T_{n,k}^{-1}$.

It remains to check that $T_{n,k}^{-1}$ carries
seeds to seeds.  If we know in advance that
$(A,B)$ is a seed, then $T_{n,k}^{-1}(A^*,B^*)$
must be a seed, namely $(A,B)$.
However, what we want to show is that $T_{n,k}^{-1}$
carries an arbitrary seed to a seed.

As readers familiar with the pentagram map know,
the inverse of the pentagram map can certainly
carry a convex polygon in the plane to a nonconvex
polygon. However, the inverse of the pentagram map
always carries a convex polygon in the projective
plane to a convex polygon in the projective
plane.  The same goes for the the seed map.

\begin{lemma}
\label{inverse}
If $(A^*,B^*)$ is an arbitrary seed of type
$(n,k)$, then the pre-image
$(A,B)=T_{n,k}^{-1}(A^*,B^*)$
is a seed of the same type.
\end{lemma}

\startproof
First we will show that $A$ is
a convex polygon in the projective plane.
$A$ is obtained from $A^*$ 
by gluing on finitely many shaded triangles,
as shown in Figure 3.5.  These
triangles are constructed successively, starting
with $A_{n-k+1}^*A_{n-k+1}^*A_{n-k+2}$
and then going counterclockwise.  A routine inductive
argument shows that these triangles are never degenerate.

We know that there are some choices of seed $(A^*,B^*)$,
namely those in the image of $T_{n,k}$, which give
rise to a strictly convex polygon $A$.  We can consider
a continuous path from a seed which has this property
to the seed we are interested in.  The fact that the
abovementioned triangles never degenerate implies
that the convexity property is both an open and
closed condition along our path of $A$-polygons.
Since the initial $A$-polygon is strictly convex,
so is the final one.

It follows from the definition of $T_{n,k}^{-1}$
that each point $B_j$ lies on the line segment
$A_jA_{j+1}$. Each edge of the form
$B_jA_j$ and $B_jA_{j+1}$ appears as an edge of
one of the shaded triangles.  Hence these
edges are all nontrivial.  This forces
$B_j$ to lie in the interior of the segment
$A_jA_{j+1}$. 
\endproof

\noindent
{\bf Remark:\/}
Projective
dualities conjugate the pentagram map (suitably
interpreted) to its inverse.  The same ought to
be true for the map $T_{n,k}$, provided that
seeds can be interpreted in a way that
puts points and lines on the same footing.
The reader who stares hard enough at Figure 3.5
will eventually see that this is possible.
Given the interaction with duality, the two
maps $T_{n,k}$ and $T^{-1}_{n,k}$ are on the
same footing and Lemma \ref{inverse} is obvious.
\newline

It we interpret ${\cal C\/}(n,k)$ as the
space of projective classes of seeds of
type $(n,k)$, then ${\cal C\/}(n,k)$ is
clearly a cell of dimension $(2n-8)+k$.
Indeed, ${\cal C\/}(n,k)$ is just a decorated
version of ${\cal C\/}(n)$.  The maps
$T_{n,k}$ and $T_{n,k}^{-1}$ both acts
as smooth diffeomorphisms on ${\cal C\/}(n,k)$.
To prove Theorem \ref{cell}, it
only remains to reconcile our definition
here with the ones made in \S \ref{maindef}.

\subsection{Reconciling the Two Definitions}

In this section we finish the proof
of Theorem \ref{cell} by reconciling the two
points of view, namely
\begin{enumerate}
\item ${\cal C\/}(n,k)$ is the space of
seeds of type $(n,k)$ and $T_{n,k}$ is the
evolution operator defined above by a
projective construction.
\item ${\cal C\/}(n,k)$ is the space of
labeled PLC pentagram spirals of type $(n,k)$
and $T_{n,k}$ is the shift map in this space.
\end{enumerate}

Suppose we start with the
first point of view. Let
$(A,B)$ be a seed.
The union of the spiral segments in the
bi-infinite orbit
\begin{equation}
\label{seedunion}
\bigcup_{q \in \Z} T_{n,k}^q(A,B)
\end{equation}
consists of $k$ embedded and locally convex bi-infinite
paths in the plane.  To see that these paths are 
pentagram spirals, we just have to see that the pentagram
map permutes them.  This is most easily seen by looking
at $3$ consecutive seeds, as in Figure 3.6.

\begin{center}
\resizebox{!}{3.5in}{\includegraphics{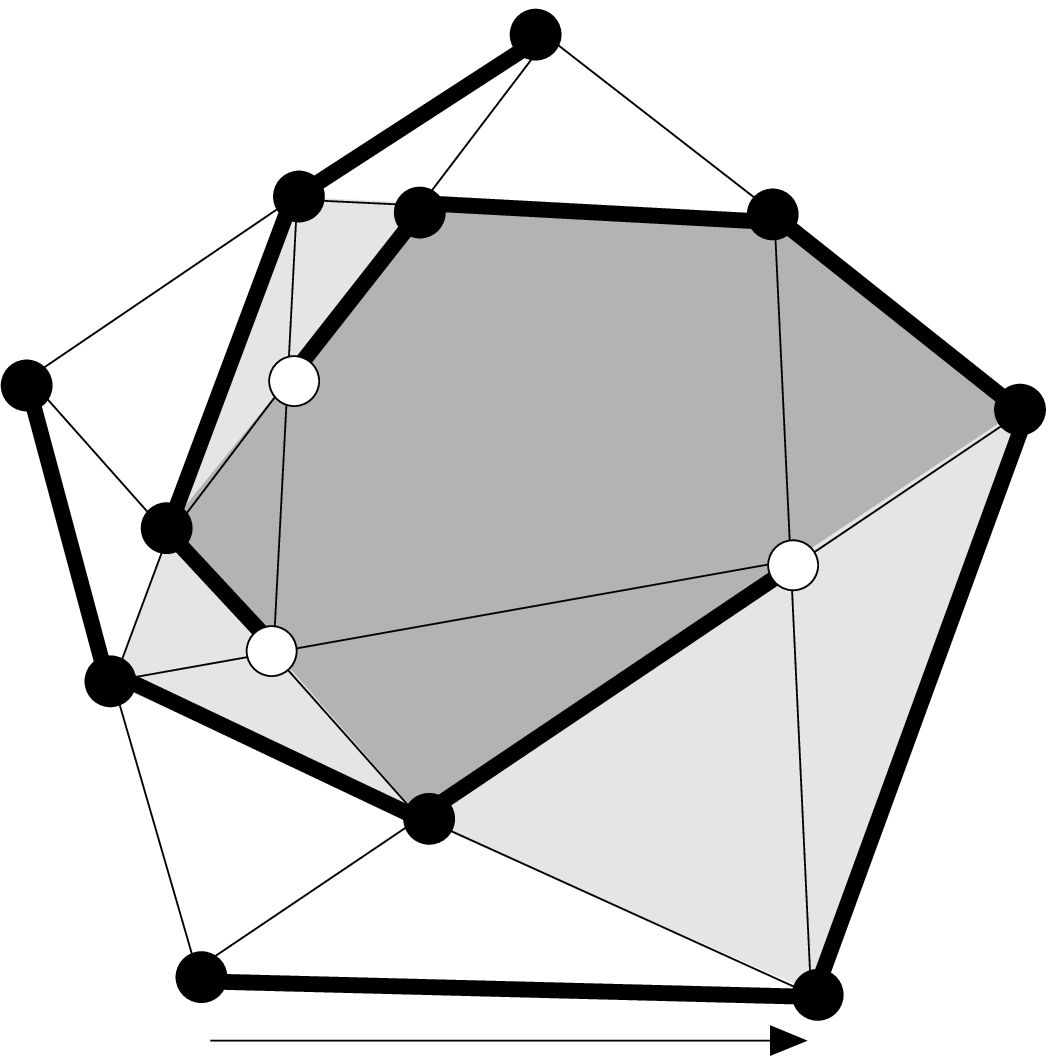}}
\newline
{\bf Figure 3.6:\/} Three consecutive seeds
\end{center}

Looking at Figure 3.6 (and generalizing)
 we see that each line drawn in
the construction of the {\it middle\/} seed is a shortest
diagonal of one of the spirals.  Since this is true for
the middle seed in any sequence of $3$ consecutive
seeds, this fact is true for every seed in the union
in Equation \ref{seedunion}.  Moreover, the white
points on each spiral are revealed to be on the
image of the ``previous'' spiral under the pentagram
map.  Since this is true for the corresponding points
of the middle seed in every consecutive run of $3$
in Equation \ref{seedunion}, we see that the pentagram
map indeed permutes the distinghished polygonal paths.
Hence, these paths are all pentagram spirals.

We label our pentagram spirals so that the distinguished
point is the point $A_1$ of the seed $(A,B)$.  We can think
of $T_{n,k}$ as acting on the union in 
Equation \ref{seedunion}.  $T_{n,k}$ takes this union
to exactly the same union, except that it is based
on the seed $(A_*,B_*)=T_{n,k}(A,B)$.  All that has happened
is that the new distinguished point is $A_1^*$, the
first point of the new seed. The arrow at the bottom
of Figure 3.6 points from $A_1$ to $A_1^*$.
With this interpretation,
$T_{n,k}$ clearly acts as the shift map on
our pentagram spirals.  

In short, every seed generates a labeled pentagram spiral,
and the seed map acts a shift on the spiral.  Thus, for
at least {\it some\/} pentagram spirals, namely those
generated by seeds, the two points of view coincide.
It remains to see that every pentagram spiral is generated
by a seed.  We will be a bit sketchy in our argument because
we don't care much about the result.  If it turned out that there
were some exotic pentagram spirals which did not come
from seeds, we would simply add the condition to our
definition that the spiral come from a seed.

\begin{lemma}
Every pentagram spiral is generated by a seed.
\end{lemma}

\startproof
Let $\Sigma$ be a pentagram spiral.  We interpret $\Sigma$ as
an infinite triangulation of some subset of the projective
plane.  The numbers $n$ and $k$ associated to $\Sigma$
characterize the global combinatorics of the tiling.
Let $\Sigma'$ be a pentagram spiral produced by a seed of
type $(n,k)$.  The tiling $\Sigma'$ has the same combinatorics
as the tiling $\Sigma$. That is, there is a homeomorphism $h$ which
carries the one triangulation to the other -- vertices are
taken to vertices and edges are taken to edges.  Moreover,
certain triples of edges in $\Sigma$ and in $\Sigma'$ line up to
form longer line segments -- the diagonals used in
the pentagram map.
$h$ respects this additional collinearity.
Evidently $h$ maps a seed for $\Sigma'$ to a seed for
$\Sigma$.
\endproof

\newpage

\section{Invariant Coordinates}

\subsection{Flags}

We mentioned in \S \ref{tiling} that we would inprove
upon the labeling scheme from \S 2.  The key is to
use flags rather than points or lines.
On a polygonal path, a {\it flag\/} is a pair $(v,e)$ where
$v$ is a vertex of the path and $e$ is an edge of the path.
We indicate the flag $(v,e)$ with an auxilliary point
placed on the edge $e$ two-thirds of the way towards $v$.
Figure 4.1 shows what we mean.

\begin{center}
\resizebox{!}{.3in}{\includegraphics{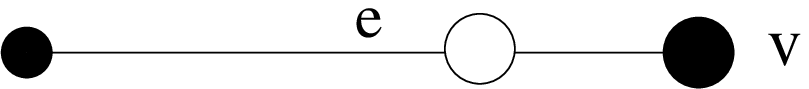}}
\newline
{\bf Figure 4.1:\/} Denoting the flag $(v,e)$.
\end{center}

Suppose we have an oriented polygonal path, as shown in 
Figure 4.2.  We orient the flags according to the
following scheme.

\begin{center}
\resizebox{!}{3in}{\includegraphics{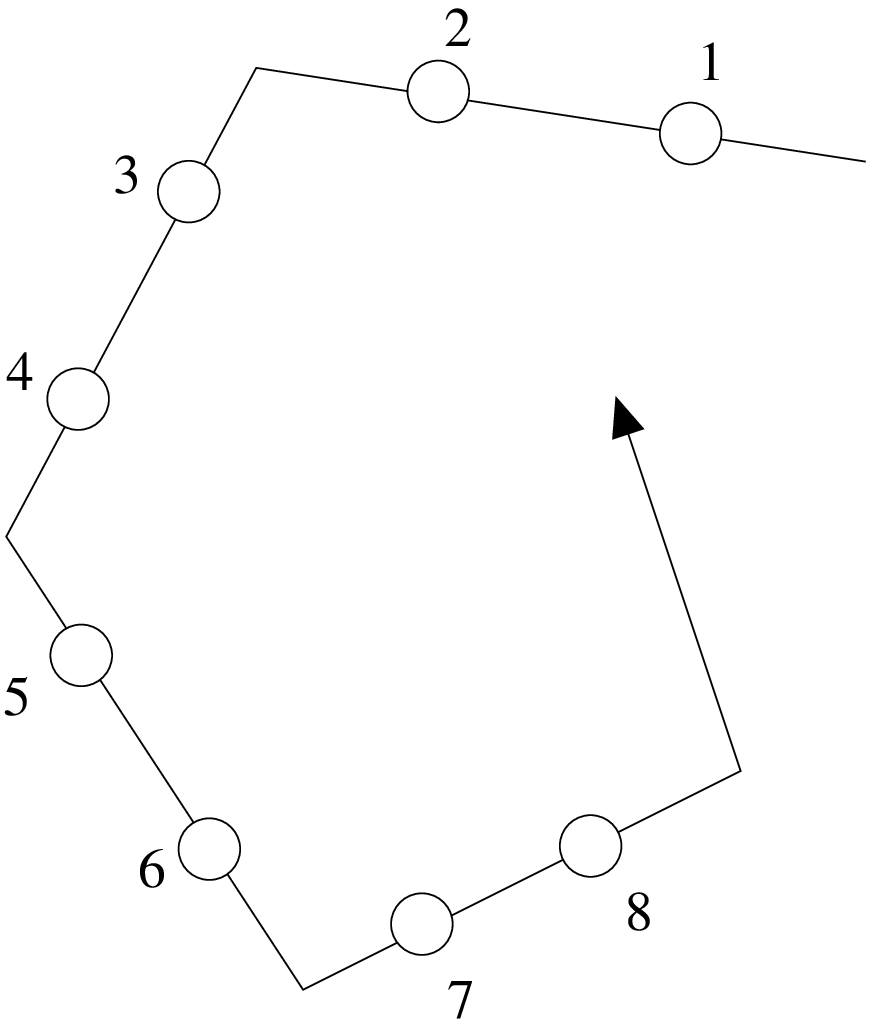}}
\newline
{\bf Figure 4.2:\/} Ordering the flags along an oriented path
\end{center}

Finally, to each flag along such a path, we associate
the cross ratio of the associated points shown in Figure 4.3.
This picture is meant to be invariant under projective
transformations.  We call these the {\it flag invariants\/}.

\begin{center}
\resizebox{!}{3in}{\includegraphics{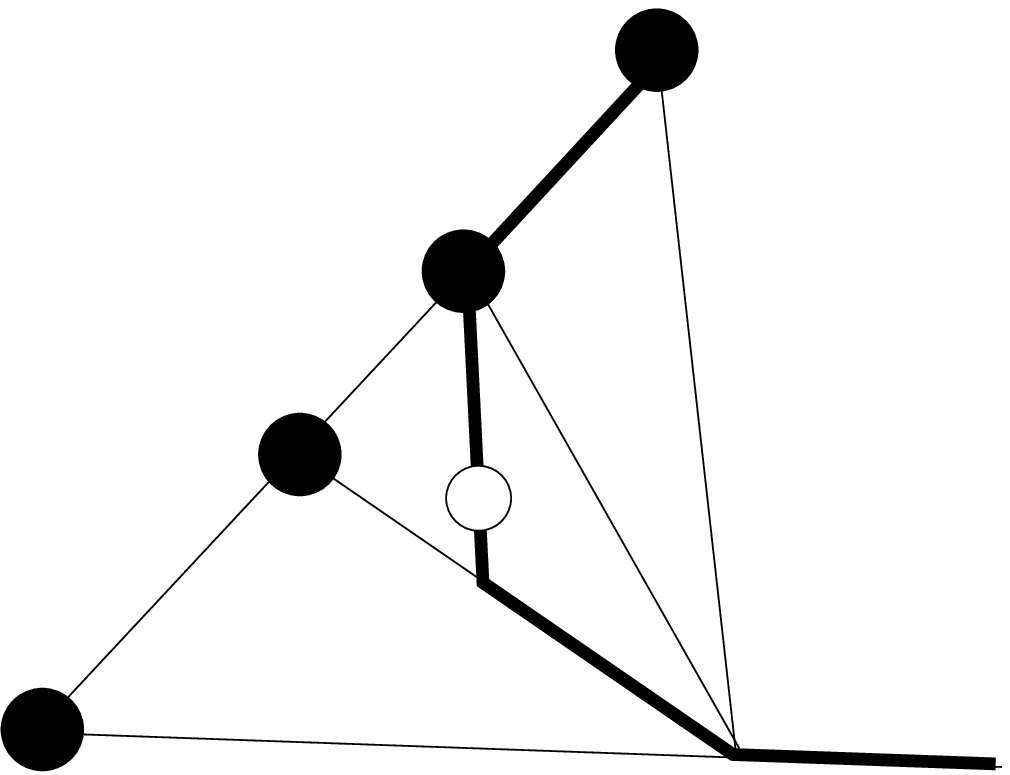}}
\newline
{\bf Figure 4.3:\/} Invariant of a flag
\end{center}

Let us consider the naturality of this construction.
The cross ratio of interest can be computed in two 
ways.  First of all, it is the cross ratio of the
$4$ points shown.  Two of the points involved are
adjacent to the flag point.  On the other hand, the
cross ratio can be computed as the cross ratio of the
$4$ drawn lines.  Two of the lines are adjacent to
the line of the flag, going in the other direction
from the abovementioned points.  Indeed, the entire
picture is invariant not just under projective transformations
but also projective dualities. Were we to apply a projective
duality to the picture, producing another polygonal curve
(with points and lines interchanged) the invariant
associated to the flag would be the same.  In short,
the invariant we have associated to the given flag is
the canonical choice.

Moreover, the ordering of the flags in Figure 4.2 reproduces
the ordering of the invariants listed in \S \ref{invariant}.
The difference is that, when we work with the flags, there
is a canonical way to line up the variables of the path
$P$ and its image $T(P)$ under the pentagram map.  Figure
4.4 below shows this.  With this new scheme, there is
a natural way to transfer the flag labels to the
hexagonal tiling so that the compatibility conditions
hold.  The method is such that same-numbered flags
correspond to diagonal edges whose centers are on
the same vertical.  This is illustrates in Figure 4.5.
With this new scheme, one need not change labeling
conventions at each level.

\begin{center}
\resizebox{!}{4in}{\includegraphics{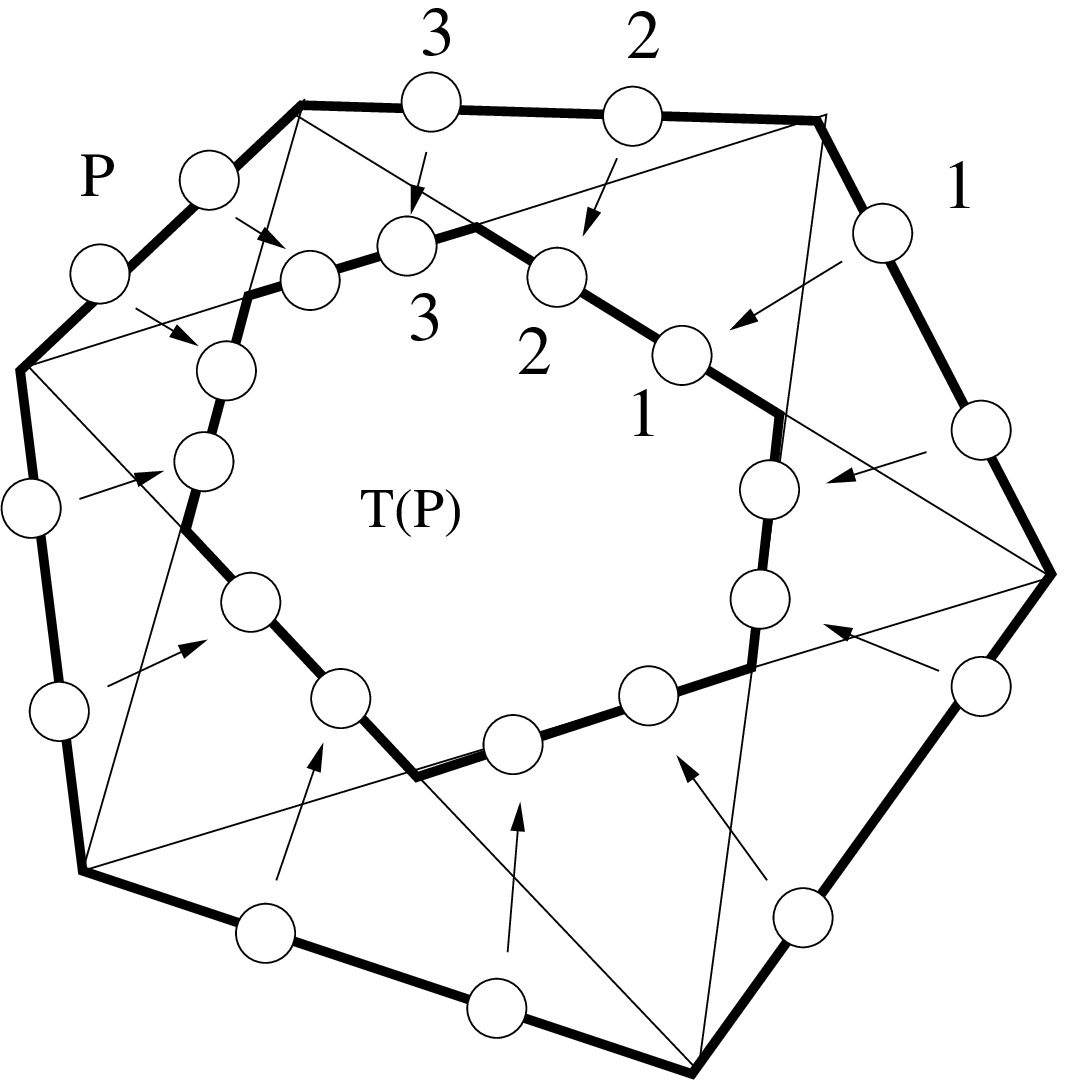}}
\newline
{\bf Figure 4.4:\/} Canonical labeling for the pentagram map.
\end{center}

\begin{center}
\resizebox{!}{2.5in}{\includegraphics{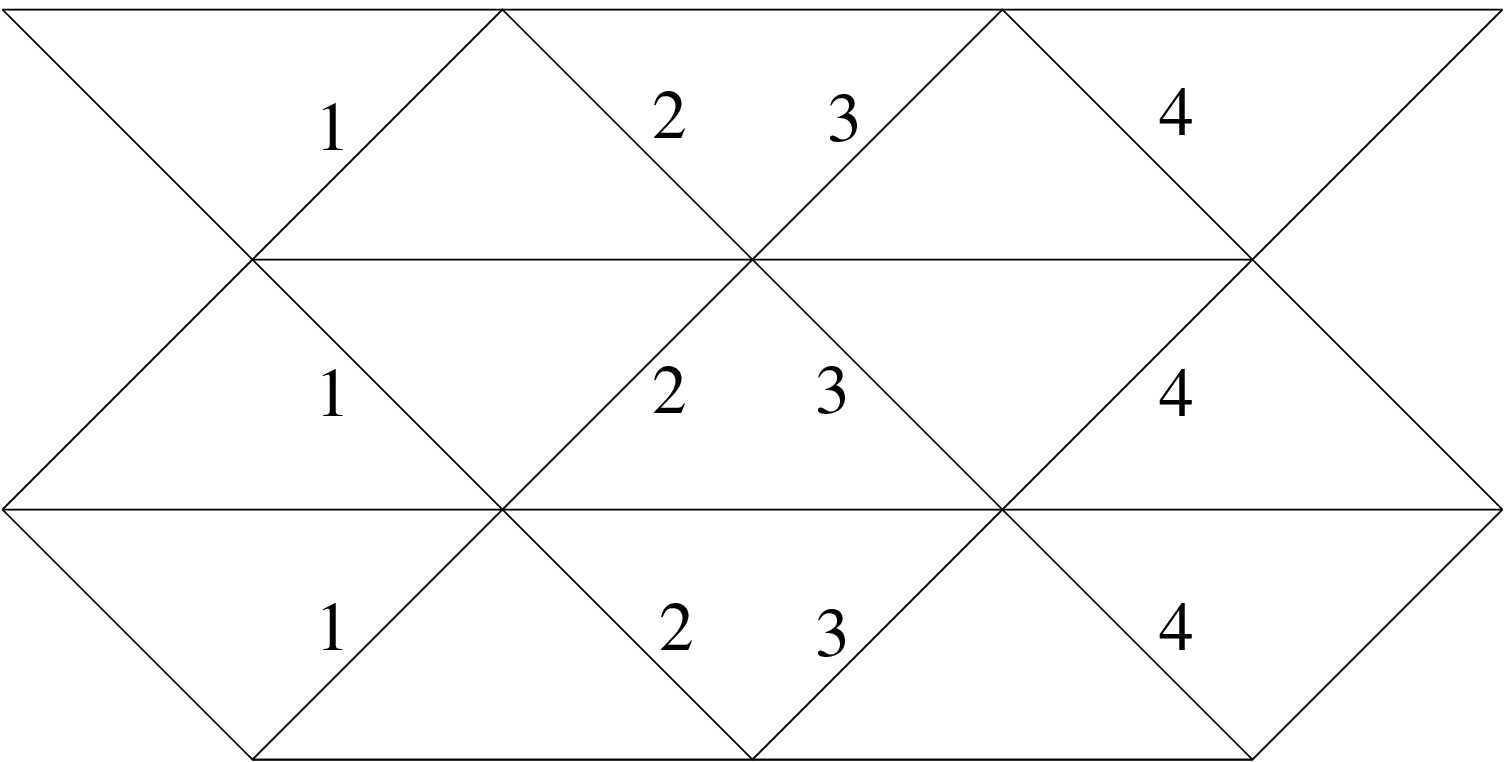}}
\newline
{\bf Figure 4.5:\/} Transfer of labels to the tiling.
\end{center}

\subsection{Pentagram Tilings Revisited}
\label{tiling2}

Let $\T$ denote the infinite triangular tiling
discussed in \S \ref{tiling}.  Again, we normalize
$\T$ so that the horizontal edges have length $1$.
The height of each triange is $1/2$.  Let
${\cal T\/}$ denote the set of all pentagram
tilings of the edges of $\T$.  

Let
${\cal T\/}(n,k)$
denote the set of pentagram tilings which are
invariant under translation by the vector
\begin{equation}
\label{trans}
V{n,k}=\Big(n+\frac{k}{2},\frac{k}{2}\Big).
\end{equation}
This vector has half-integer coordinates,
and the sum of the two coordinates is an integer.
Hence, translation by $V_{n,k}$ is an isometry of $\T$.
Thus, the definition of
${\cal T\/}(n,k)$ is not vacuous.

Comparing the discussion in \S \ref{maindef} with the
scheme for transfering the flag invariants of
flags to the labels of $\T$, we see that
each element of ${\cal C\/}(n,k)$ gives rise
to a unique element of ${\cal T\/}(n,k)$. 
To make this completely precise, we translate
$\T$ so that one vertex lies at the origin.
We arrange that the diagonal edge joining
$(0,0)$ to $(1/2,1/2)$ (respectively to
$(-1/2,-1/2)$)
corresponds to the flag just
inward (respectively outward)
from the distinguished vertex of the
spiral.  Thus, when we move righward
along horizontal edges, it corresponds to
going inward along the spiral.
When we move downward, it corresponds to
doing the pentagram map.

It is not true that every element of
${\cal T\/}(n,k)$ arises from an element
of ${\cal C\/}(n,k)$.  For one thing,
all the labels would have to lie in $(0,1)$.
See \S \ref{convexity}.
For another thing, the dimensions of the
spaces do not match up.  For instance, in the
toy case when $k=0$, the former space has
dimension $2n$ and the latter space has
dimension $2n-8$.  The other elements in
${\cal T\/}(n,0)$ correspond to the so-called
{\it twisted polygons\/}.  In \S 7 we will
informally discuss a similar 
interpretation of the general
element of ${\cal T\/}(n,k)$. Here we will only
consider elements of
${\cal T\/}(n,k)$ which come from elements
of ${\cal C\/}(n,k)$.
\newline
\newline
{\bf Remark:\/}
The elements of ${\cal T\/}(n,k)$ can be considered as
edge labelings of a triangulation of the
cylinder $\R^2/V_{n,k}$.
The combinatorics 
of this triangulation is essentially the
same as the combinatorics of the tiling
obtained from the pentagram spiral, though
the edges of the one triangulation do not
precisely match up with the edges of the other.
We leave it to the interested reader to work
out the exact correspondence. 

\subsection{Discussion of Formulas}

It might be nice, or at least useful for further research,
to give explicit formulas for the action of $T_{n,k}$ on
${\cal C\/}(n,k)$.  This amounts to
identifying ${\cal C\/}(n,k)$ with a specific algebraic
variety, and then expressing $T_{n,k}$ as a birational
transformation of that variety.
Essentially all the papers on the pentagram 
map take this approach.

We have seen above that each element of ${\cal C\/}(n,k)$
gives rise to an edge-labeling of a
triangulation of $\R^2/V_{n,k}$.
The labeling satisfies the above
compatibility rules,
and there is some finite list ${\cal L\/}$ of
edges such that the labels on $\cal L$ determine
all the other labels.  Thus, one can realize
${\cal C\/}(n,k)$ as an algebraic variety
in a finite dimensional space.
Using the compatibility rules, one can 
express $T_{n,k}$ as a birational map.

For the pentagram map, this approach is
completely successful. There is a canonical
choice for $\cal L$, and the compatibility rules
give rise to a transformation with a very nice
formula.  See [{\bf Sch3\/}], [{\bf OST1\/}],
and [{\bf OST2\/}].  However, for the pentagram
spirals, I have not been able to find a 
good choice for $\cal L$.  No choice seems
canonical, and all choices seem to 
lead to messy formulas. 

This state of affairs does not (yet) bother me.
I think that the right point of view is that
elements of ${\cal C\/}(n,k)$ are simply these
labeled triangulations with the compatibility
rules and the map $T_{n,k}$ is just a shift
operator on the space of such labelings.
However, I can see that this answer will be
unsatisfying to some readers, and perhaps I
will eventually find it unsatisfying.  
I hope with the interested reader will take
up the question of finding good formulas in
the sense discussed above.  This main point of
the discussion is that the problem is nontrivial.

\subsection{The First Invariant}
\label{first}

We will use a pictorial 
method for representing polynomial invariants
of the map $T_{n,k}$.  Let $C$ be some
collection of edges of $\T$.  The product of
the labels of the variables associated to the
edges constituting $C$ is a monomial which we denote by
$\langle C \rangle$.  So, $\langle C \rangle$ denotes a
function on ${\cal T\/}(n,k)$ defined by $C$.

Given a vector $V$, we say that
$\langle C \rangle$ is $V$-{\it invariant\/} if
the two functions
$\langle C \rangle$ and
$\langle C+V \rangle$
agree on ${\cal T\/}(n,k)$.
Here $C+V$ is the copy of $C$ which has
been translated by $V$.  This definition
only makes sense when translation by $V$
preserves $\T$. We call such vectors
{\it allowable\/}.
Though $n$ and $k$ are not explicitly
mentioned, it is understood that the
notion of invariants is defined with
respect to these parameters.

Figure 4.5 shows the configuration $C$ 
associated to the pair $(n,k)=(4,3)$.
This configuration ``goes up'' $4+3=7$ steps
and then ''does down'' $4$ steps.
Translation by $V_{4,3}$ identifies the
endpoints of $C$ and thus $C$ defines a
closed path on the cylinder $\R^2/V_{4,3}$ discussed in the
remark at the end of the last section.
We will use the notation $Z(4,3)$ to
denote the corresponding function $\langle C\rangle$
in this case.  The general definition of $Z(n,k)$ 
follows the same pattern.

\begin{center}
\resizebox{!}{2.3in}{\includegraphics{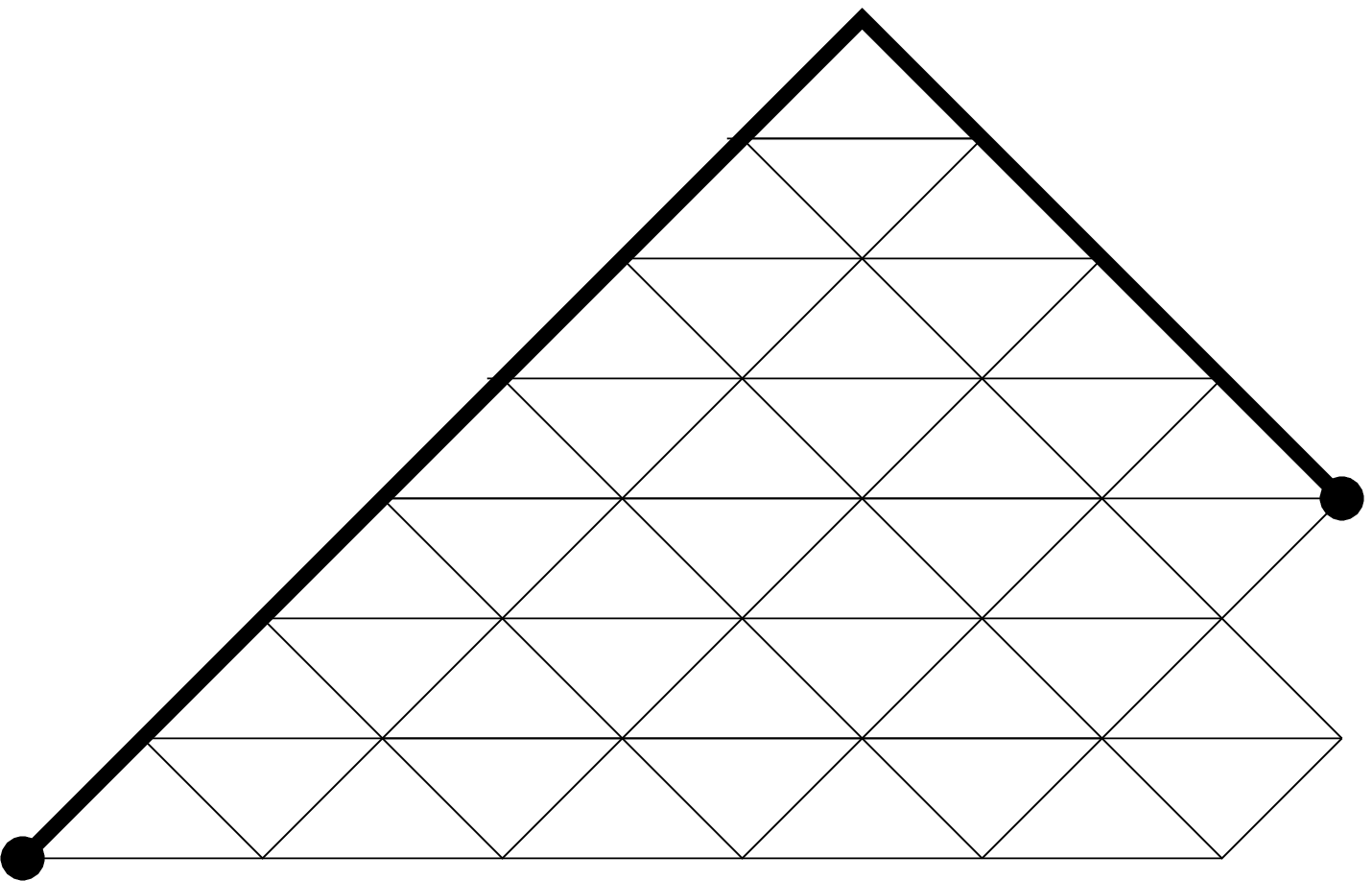}}
\newline
{\bf Figure 4.5:\/} The monomial $Z(4,3)$.
\end{center}

Below we will prove that $Z(n,k)$ is invariant
with respect to any vector of $\T$.  Before we prove
this result, we need to make a short digression.
Say that a {\it zigzag\/} is a path which
moves rightward, along diagonal edges of
$\T$ and joins two vertices of $\T$.

\begin{lemma}
Suppose that $Z_1$ and $Z_2$ are zigzags which 
start and end at the same vertex.  Then
$\langle Z_1 \rangle = \langle Z_2 \rangle$.
\end{lemma}

\startproof
If we push a
zigzag across a single square, as shown in Figure 4.6
below, the corresponding monomial does not change,
thanks to the pentagram relations.  So, we just
keep pushing the one zig-zag until it equals the other.
\endproof

\begin{center}
\resizebox{!}{1in}{\includegraphics{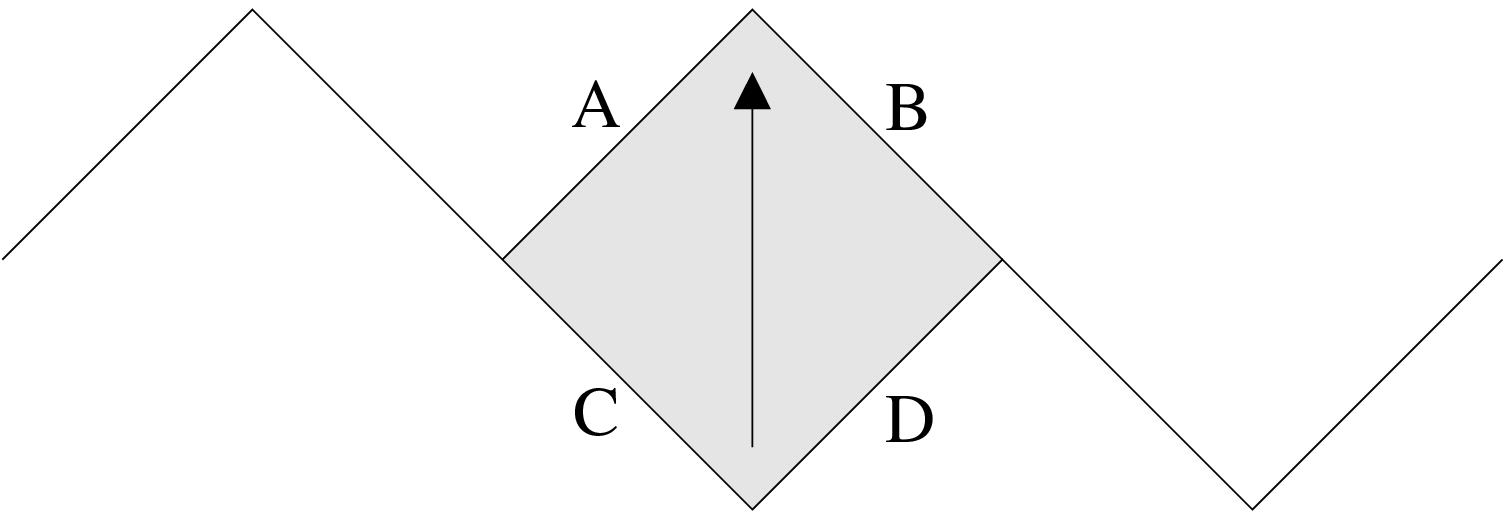}}
\newline
{\bf Figure 4.6:\/} pushing a zigzag: $AB=CD$.
\end{center}

\begin{lemma}
$Z(n,k)$ is invariant with respect to the vector
$(1,0)$.
\end{lemma}

\startproof
Let $C$ be such that
$Z(n,k)=\langle C \rangle$ and let $C'=C+(1,0)$.
Let $Z'(n,k)=\langle C' \rangle$.
Referring to Figure 4.7, we have
\begin{equation}
C=A \cup B; \hskip 30 pt
C'=A' \cup B'.
\end{equation}

\begin{center}
\resizebox{!}{2.6in}{\includegraphics{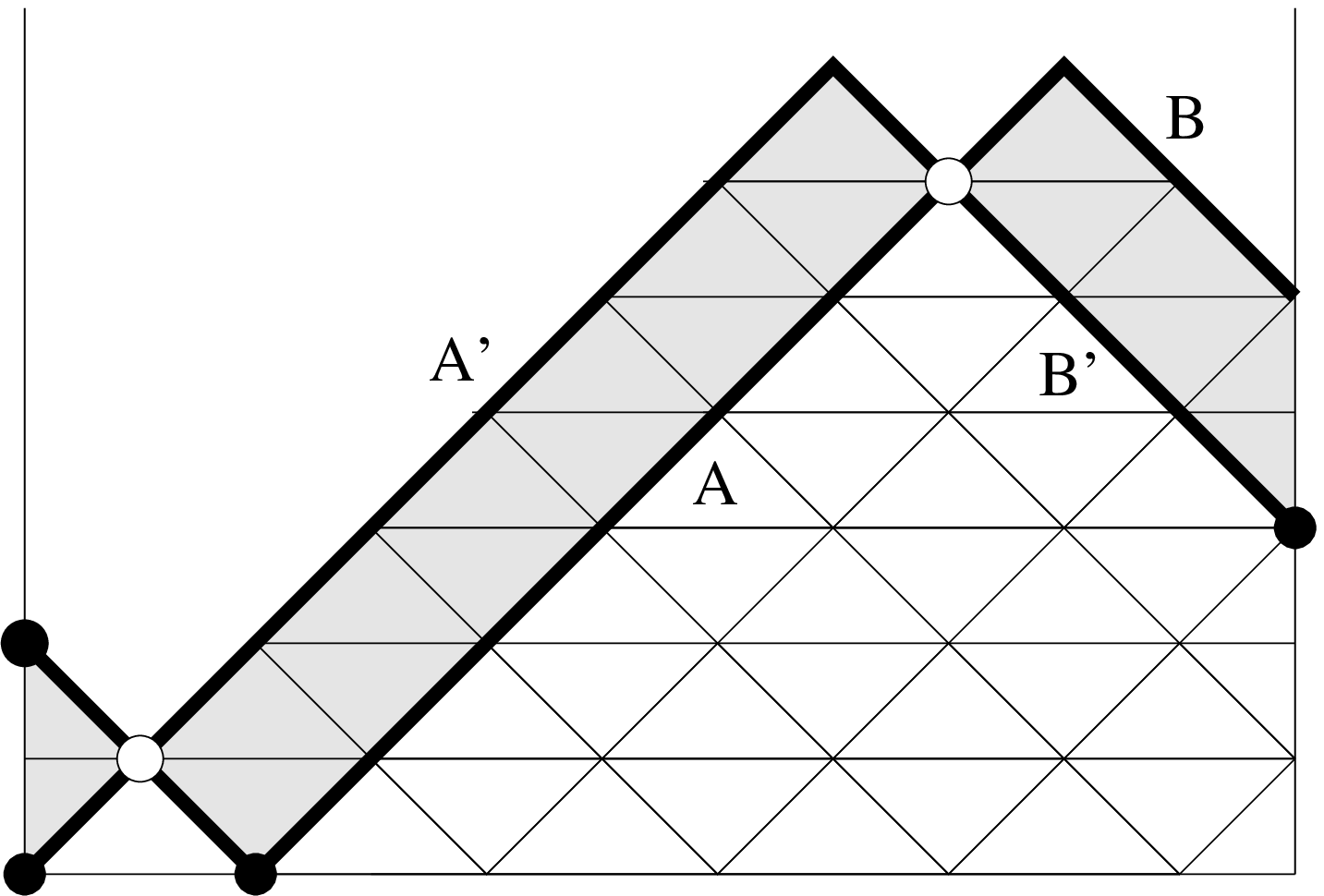}}
\newline
{\bf Figure 4.7:\/} $C \cup C'$ drawn on the cylinder
$\R^2/V_{4,3}$.
\end{center}

By the lemma, we have
$\langle A \rangle = \langle A' \rangle$ and
$\langle B \rangle = \langle B' \rangle$. Finally,
we have
$Z(n,k)=\langle A \rangle \langle B \rangle$ and
$Z'(n,k)=\langle A' \rangle \langle B' \rangle$.
\endproof

Essentially the same argument works for
the vector $(1/2,1/2)$.  Since the
vectors $(1,0)$ and $(1/2,1/2)$ generate
$\T$, we see that $Z(n,k)$ is invariant
with respect to any allowable $V$.

\newpage

\section{Compactness of the Orbit Closures}
\subsection{Local Convexity}
\label{convexity}

Our main goal in this chapter is to prove
Theorem \ref{compact}.  For ease of exposition
we will assume that $(n,k) \not = (4,1)$.
Theorem \ref{period} from the next chapter
takes care of this exceptional case.

\begin{lemma}
\label{pentagon}
Let $P$ be a PLC pentagram spiral whose type
is not $(4,1)$.
Every $5$ consecutive points of $P$ are the vertices
of a strictly convex pentagon.
\end{lemma}

\startproof
Let $(n,k)$ be the type of $P$.
It suffices to consider the
points $P_1,P_2,P_3,P_4,P_5$ and let $(A,B)$ be the
seed such that $A_1=P_1$ and $A_2=P_2$.   There are
several cases.
 
If $n-k \geq 4$, then
$P_i=A_i$ for all $i$, and the result is clear:  $A$ is a
strictly convex polygon.  

If $n-k=3$ (and $n>4$) then
$P_i=A_i$ for $i=1,2,3,4$ and $P_5=B_4$. Again, by
construction, the result is true.

\begin{center}
\resizebox{!}{2.7in}{\includegraphics{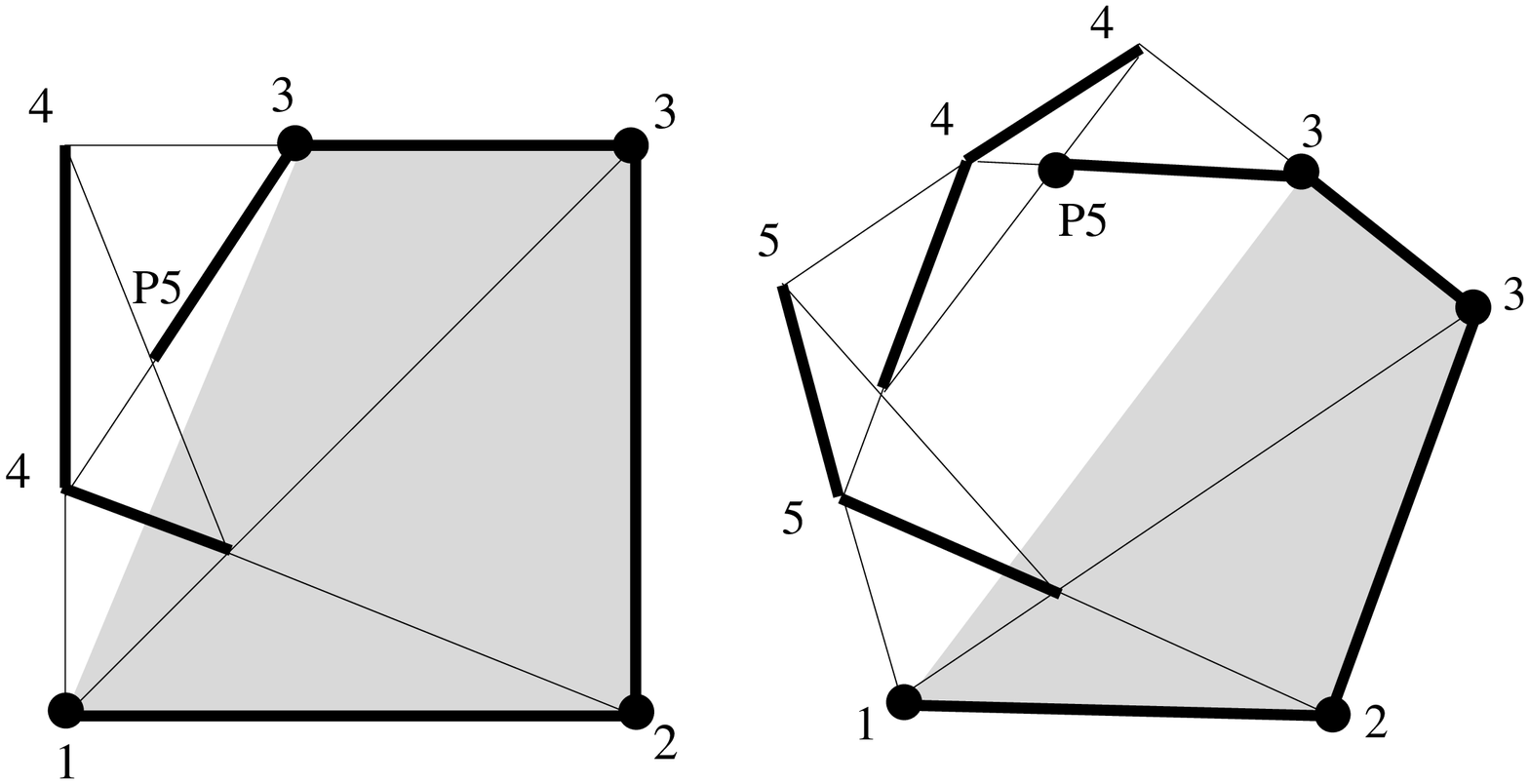}}
\newline
{\bf Figure 5.1:\/} The case when $n-k=2$ for $n=4,5$.
\end{center}

If $n-k=2$ then
$P_i=A_i$ for $i=1,2,3$ and $P_4=B_3$.  The points
$P_1,P_2,P_3,P_4$ therefore
form the vertices of a convex quadrilateral $Q$, shaded
in Figure 5.1.
The line $P_4P_5$, which
coincides with the line $B_3B_4$, lies
outside $Q$ and inside
the convex polygon bounded by $A$.
These two properties imply our result.

\begin{center}
\resizebox{!}{2.7in}{\includegraphics{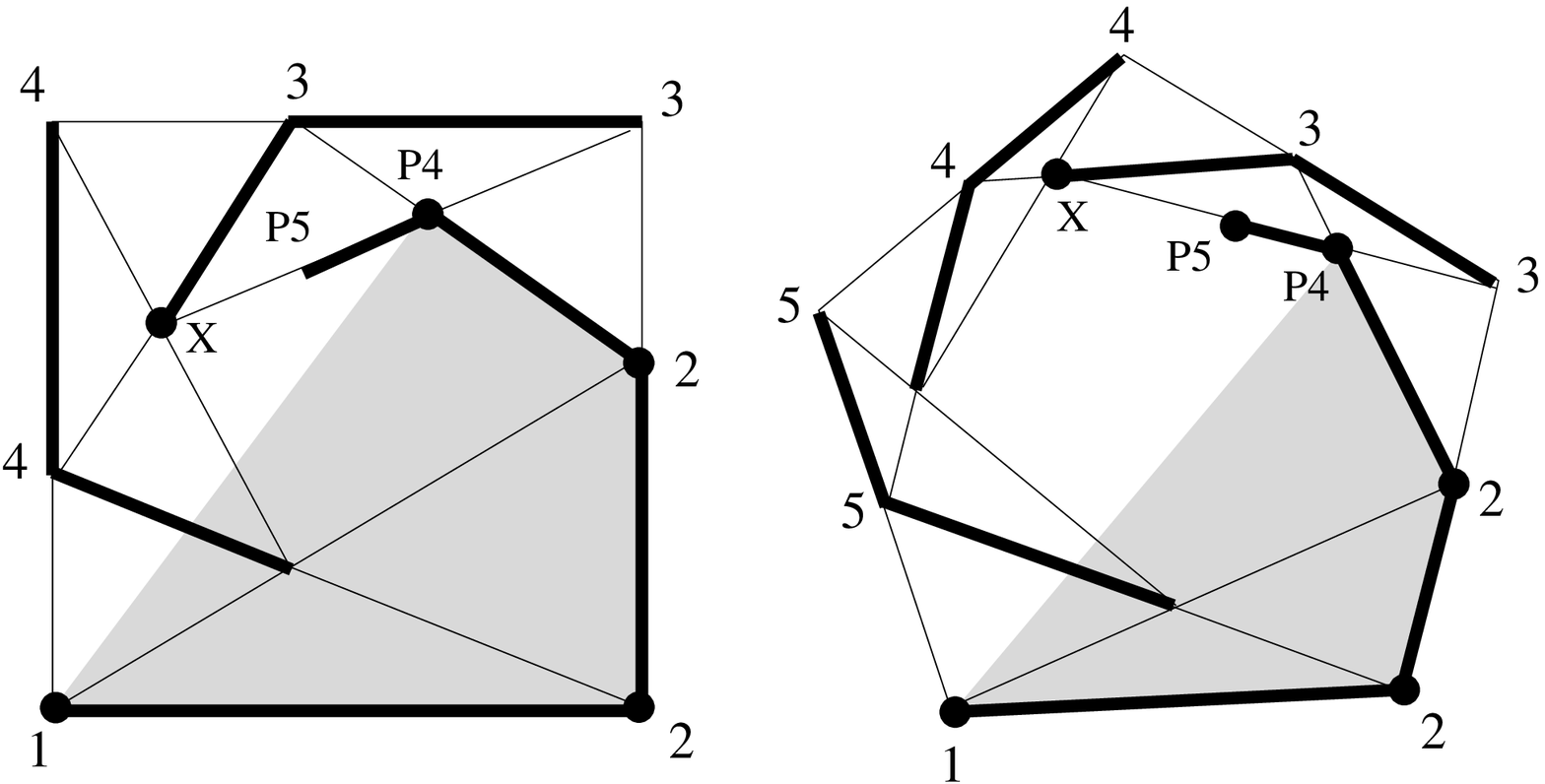}}
\newline
{\bf Figure 5.2:\/} The case when $n-k=1$ for $n=4,5$.
\end{center}

If $n-k=1$ then
$P_i=A_i$ for $i=1,2$ and $P_3=B_2$.  The same argument as
in the preceding section shows that
$P_1,P_2,P_3,P_4$ form the vertices of a convex quadrilateral
$Q$, shaded in Figure 5.2.  The segment $B_3B_4$, which
contains the point marked $X$, lies outside $Q$.  Hence
$X$ lies outside $Q$ as well.  But then $P_5$, which
lies on the line segment $P_4X$, lies outside $Q$ as
well. Finally, the segment $P_4A$ lies inside the
convex polygon bounded by $A$.  These properties imply
the result.
\endproof

Given a pentagram spiral $P$ of type $(n,k)$, let
$Z(P)$ denote the value of the invariant
$Z(n,k)$ evaluated on $P$.

\begin{corollary}
\label{lower}
Suppose that $P$ is a PLC pentagram spiral.
The flag invariants associated to $P$
all lie in $(0,1)$. Hence
$Z(P)$ serves as a lower bound for 
the flag invariants.
\end{corollary}

\startproof
Each flag invariant is computed using $5$ points, as in
Figure 4.3. The convexity of the $5$ points guarantees
that the $4$ points relevant for the cross ratio come
in order along the line.  Hence, all the flag invariants
lie in $(0,1)$.  Given the invariance of $Z(n,k)$ with
respect to any allowable vector, as established
in \S \ref{first}, we see that $Z(P)$ can be
computed as a product of flag invariants, one of 
which is the one that currently interests us.
\endproof

\subsection{The Vertex Invariant}

In this section we associate a different invariant
to the the vertices (as opposed to the flags)
of a locally convex polygonal path.

\begin{center}
\resizebox{!}{2.8in}{\includegraphics{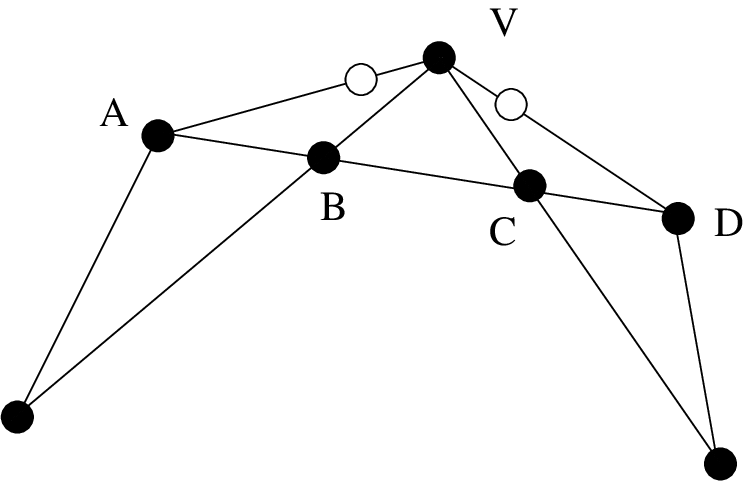}}
\newline
{\bf Figure 5.3:\/} The vertex invariant
\end{center}

Referring to the points in Figure 5.3, the
{\it vertex invariant\/} is
\begin{equation}
\chi(v)=[A,B,C,D].
\end{equation}
A routine calculation,
which we omit, shows that
\begin{equation}
\label{chi}
\chi(v)=f_1f_2,
\end{equation}
where $f_1$ and $f_2$ are the invariants associated
to the flags adjacent to $v$.  Such a relation
is not so surprising, because everything in sight
just depends on the $5$ points shown.

\begin{lemma}
\label{bound}
Let $P$ be a PLC pentagram spiral.
The quantity $Z^2(P)$ is a lower bound
for the vertex invariants of $P$.
\end{lemma}

\startproof
This is immediate from 
Corollary \ref{lower} and from 
Equation \ref{chi}.
\endproof

\subsection{Uniform Bounds}

One of the main goals in this chapter is to show that
the invariant $Z$ has compact level sets.  In this
section, we will consider a sequence
$\{P(m)\}$ of pentagram spirals such that $Z(P(m))$ is
independent of $m$.  Our goal is to establish some uniform
bounds for such spirals.  We will normalize by projective
transformations so that 
\begin{equation}
P_1(m)=(0,0), \hskip 15 pt
P_2(m)=(1,0), \hskip 15 pt
P_3(m)=(1,1), \hskip 15 pt
P_4(m)=(0,1)
\end{equation}
for all $m$.  Here $P_j(m)$ denotes the $j$th point of $P(m)$.
We let $P_+(m)$ denote the union of all points $P_j(m)$ where
$j \geq 1$.  This is an inward spiraling half of $P(m)$.

\begin{lemma}
\label{long}
There is some $D$ such that
$\|P_5(m)-P_4(m)\|<D$.  Here
$D$ does not depend on $m$.
\end{lemma}

\startproof
It follows from Lemma \ref{pentagon} and from our
normalization that $P_5(m)$ lies between the two lines
$y=0$ and $y=1$, as shown in Figure 5.3.

\begin{center}
\resizebox{!}{2.8in}{\includegraphics{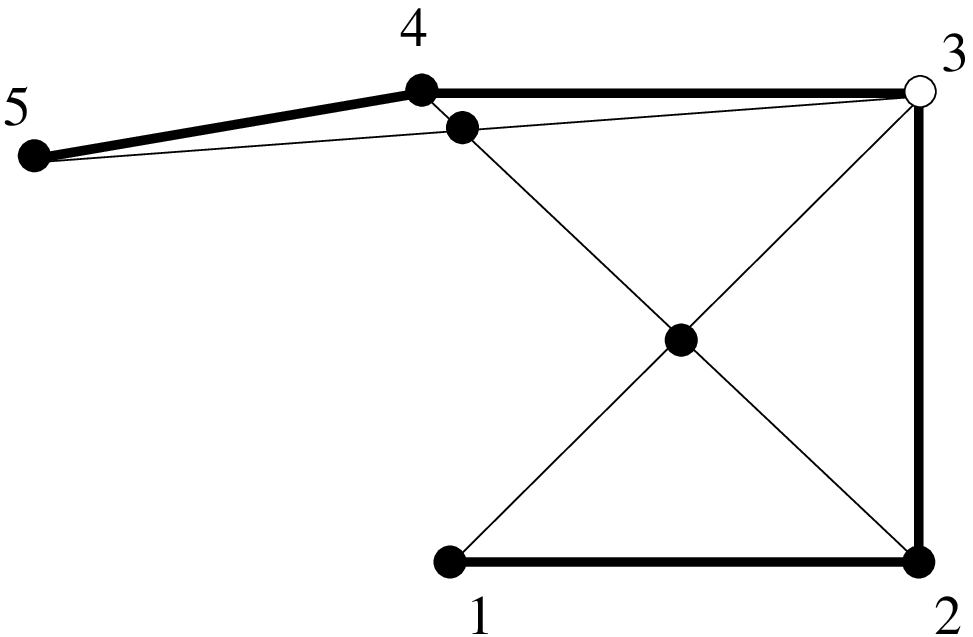}}
\newline
{\bf Figure 5.3:\/} The points $P_1,P_2,P_3,P_4,P_5$.
\end{center}

If
$\|P_5(m)-P_4(m)\| \to \infty$ then the corner invariant
$\chi(P_3)$ tends to $0$ with $m$.  This is impossible,
by Lemma \ref{bound}.
\endproof

\begin{lemma}
There is some $d>0$ such that
$\|P_5(m)-P_4(m)\|>d$. Here
$d$ does not depend on $m$.
Likewise, there is some $s>0$ such that
the line $P_4(m)P_5(m)$ has slope at
least $s$.  Here $s$ is independent of $m$.
\end{lemma}

\startproof
These statements have the same proof as in Lemma
\ref{long}.
\endproof

\begin{lemma}
There is some compact subset $K \subset \R^2$ such that
$P_+(m) \subset K$ for all $m$.  In particular, there
is some $D$ such that $\|P_j(m)-P_{j+1}(m)\|<D$ for
all $m$ and all $j>0$.
\end{lemma}

\startproof
The local convexity of $P(m)$ guarantees that
$P_+(m)$ is contained in the quadrilateral
$K(m)$ bounded by the lines
$x=1$ and $y=0$ and $y=1$ and $P_4(m)P_5(m)$.
Given what we have established about 
$P_5(m)$, we have a uniform upper bound to the
diameter of this quadrilateral.  Hence there is
a single compact $K$ which contains $K(m)$
for all $m$.  The second conclusion of the
lemma is now immediate.
\endproof

\begin{lemma}
Suppose that $j>1$.
Then $\|P_{j+1}(m)-P_j(m)\|>d$. Here
$d$ does not depend on $m$.
Likewise, the segment $P_j(m)P_{j+1}(m)$ makes
an exterior angle of at least $\theta>0$
with $P_{j-1}(m)P_j(m)$. Again, $\theta$ does
not depend on $m$.
\end{lemma}

\startproof
If this result fails, there is some smallest $j$ where the
problem goes wrong.   This means that the points
$P_{j-3},P_{j-2},P_{j-1},P_j$ are spaced uniformly far
apart independent of $m$.   Moreover, there is a uniform
lower bound to the exterior angles shown in Figure 5.4.
But then the same argument as in the proof of Lemma \ref{long},
applied to $\chi(P_{j-1})$,
gives a contradiction.
\endproof

\begin{center}
\resizebox{!}{1.8in}{\includegraphics{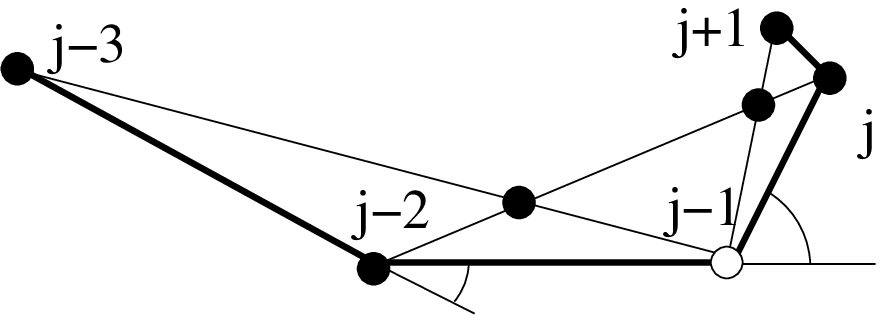}}
\newline
{\bf Figure 5.4:\/} The points $P_{j-3},...,P_{j+1}$.
\end{center}

\begin{lemma}
Suppose that $j>1$.
The segment $P_j(m)P_{j+1}(m)$ makes
an interior angle of at least $\theta>0$
with $P_{j-1}(m)P_j(m)$. Here, $\theta$ does
not depend on $m$.
\end{lemma}

\startproof
This is forced by the fact that every $5$ consecutive points
lie on a convex pentagon, together with the uniform lower
bound on the side lengths.
\endproof

\begin{corollary}
\label{converge}
The sequence $\{P_+(m)\}$ converges, at least
 on a subsequence,
to a strictly convex infinite polygonal path $P_+(\infty)$.
Every $5$ consecutive points of $P_+(\infty)$ are the
vertices of a strictly convex pentagon.  
\end{corollary}

\startproof
This follows immediately from the uniform bounds
we have on all the lengths and angles, together
with compactness.
\endproof

\subsection{Compactness Proof}

Now we put everything together and prove Theorem \ref{compact}.
As we mentioned above, it suffices to prove that the
invariant $Z=Z(n,k)$ has compact level sets.
Let $\{P(m)\}$ be as in the previous section.
These pentagram spirals all have the same $Z$-invariant.
Let $X_+(m)$ denote the portion of the triangulation
associated to $P(m)$ which is contained in the seed
\begin{equation}
\label{seed}
S'(m)=T^{4k+4}_{n,k}(S(m)).
\end{equation}
We use the seed $S'(m)$ rather than $S(m)$ because
our construction produces a union of triangles whose
outermost portion is somewhat ragged.

We will show below that the sequence $\{X_+(m)\}$ converges
to a triangulation which we will call
$X_+(\infty)$.  There are two kinds of edges in
 $X_+(m)$.  The edges of the first kind are contained
in $P(m)$ and its images under powers of the pentagram  map.
The union of these edges
will converge to strictly locally convex paths in
$X_+(\infty)$.  The second kind of edge in $X_+(m)$ is
a shortest diagonal of $P(m)$ or one of its
images under powers of the pentagram map.  These
edges of the second kind will
converge to the corresponding shortest diagonals of
the strictly convex paths in $X_+(\infty)$
which we have just mentioned.

The convergence of triangulations implies that the sequence of seeds
$\{S'(m)\}$ converges to a nontrivial seed $S'(\infty)$. 
Indeed, the convex polygon $A'$ supporting $S'(\infty)$ is
just the outer boundary of the tiling $X_+(\infty)$.
Some of the edges of this outer boundary are edges
of the first kind, and these determine the $B'$ points
marking $k$ of the edges of $A'$.  The convergence
of seeds implies the compactness result.
\newline

For ease of exposition we suppose $k>1$. The
case $k=1$ is similar and in fact easier.
Each pentagram spiral $P(m)$ is part of a
system of $k$ pentagram spirals $P(m,j)$ for $j=1,...,k$,
which are permuted by the pentagram map. 

\begin{center}
\resizebox{!}{1.5in}{\includegraphics{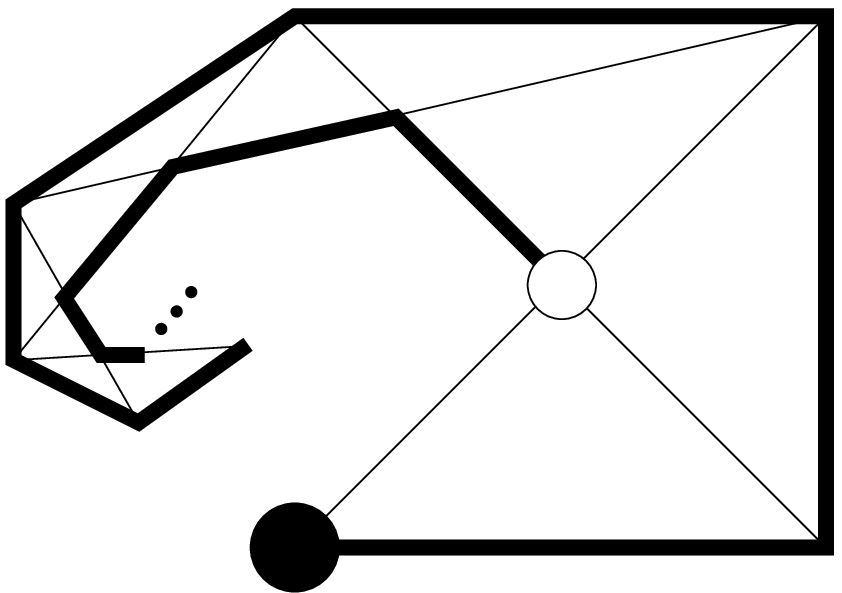}}
\newline
{\bf Figure 5.5:\/} The image of $P_+(\infty)$ under
the pentagram map.
\end{center}

Let $P_+(\infty,1)$ be the limit guaranteed by
Corollary \ref{converge}.
The pentagram map is well defined on $P_+(\infty,1)$,
starting from the first point onwards, as indicated
by Figure 5.5. 
Let $P_+(\infty,2)$ be the image of $P_+(\infty,1)$ under
the pentagram map.
Evidently, $P_+(\infty,2)$ is the limit of 
$\{P_+(m,2)\}$, where $P_+(m,2)$ is some
forward-infinite portion of $P(m,2)$.
From our construction, we see that
$P_+(\infty,2)$ is a locally strictly convex infinite path.
Moreover, every $5$ consecutive points of $P_+(\infty,2)$
are the vertices of a strictly convex pentagon.
But now the pentagram map is defined on
$P_+(\infty,2)$, and we arrive at $P_+(\infty,3)$, which
is evidently the limit of a suitably defined
sequence $\{P_+(m,3)\}$.

We continue this way until we reach
$P_+(\infty,k+1)$.  This is a proper sub-path
of $P_+(\infty,1)$, where roughly $n-k/2$ vertices
have been chopped off the front.
Let $X'_+(\infty)$ denote the union of the paths
$P_+(\infty,j)$, for $j=1,...,k$, together with
all their shortest diagonals.  
The corresponding union $X'_+(m)$ contains
the tiling $X_+(m)$. 

By construction, $X_+'(m)$ converges to
$X_+'(\infty)$.  But then there is a subset
$X_+(\infty)$ which is the limit of the
slightly smaller $X_+(m)$.  The outer boundary
of $X_+(\infty)$ is the limit of the sequence
$\{S'(m)\}$ of seeds.  This convergence is
what we had wanted to establish.  The proof
of Theorem \ref{compact} is done.

\newpage
\section{Geometry of the Tiling}

\subsection{Hilbert Diameter}

Our goal in this chapter is to prove
Theorems \ref{omega} and \ref{wind}.

Let $K \subset \R\P^2$ be a compact
convex domain.  Given $2$ points
$b,c \in K$, we have the Hilbert distance
\begin{equation}
d_K(b,c)=-\log [a,b,c,d]
\end{equation}
where $a$ and $d$ are the two points where
the line $bc$ intersects $\partial K$,
ordered as in Figure 6.1.
This is a projectively natural metric on $K$.
When $K$ is a circle, $d_K$ is the hyperbolic
metric in the Klein model.

\begin{center}
\resizebox{!}{2in}{\includegraphics{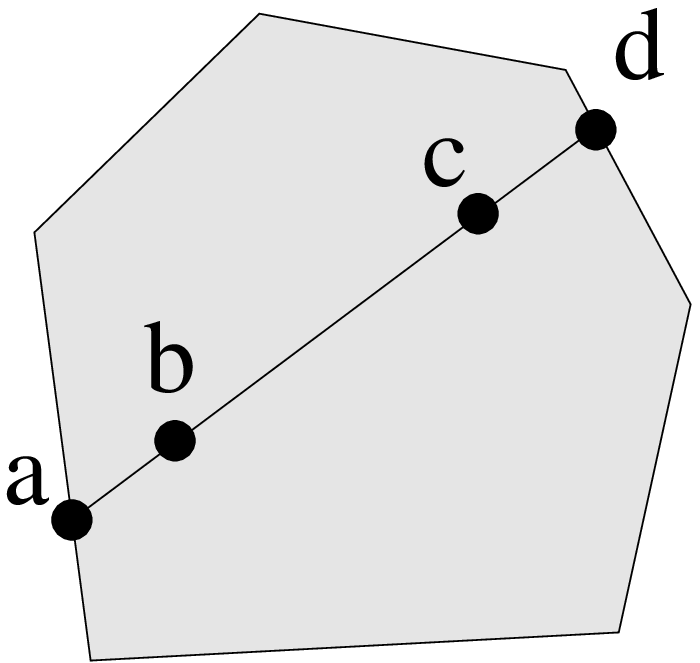}}
\newline
{\bf Figure 6.1:\/} The Hilbert distance
\end{center}

Suppose that $L$ is a compact set contained in the
interior of $K$.  We define the {\it Hilbert diameter\/}
of $L$ to be the diameter of $L$ as measured in the
Hilbert metric. 

\begin{lemma}
\label{nest}
Suppose that $\{L_n\}$ is a sequence of compact
convex subsets contained in the interior of $K$.  Suppose
that the Euclidean diameter of $L_n$ converges to
the Euclidean diameter of $K$.  Then the
Hilbert diameter of $L_n$ converges to $\infty$.
\end{lemma}

\startproof
Take a segment $\sigma_m$ which connects
two points on $L_m$ which are maximally spaced apart.
Call these points $b_m$ and $c_m$.  Let $a_m$ and
$d_m$ be the points used in the definition of
$d_K(a_m,b_m)$.  By construction
$\|a_m-b_m\|$ and $\|c_m-d_m\|$ converge to $0$ whereas
$\|b_m-c_m\|$ is uniformly bounded away from $0$.
In this situation $d_K(b_m,c_m) \to \infty$.
\endproof
\subsection{Support of the Tiling}

Now we turn to the proof of Theorem \ref{omega}.
We begin with a corollary of Lemma \ref{nest}.

\begin{corollary}
\label{shrink}
Suppose that $\{K_m\}$ is a nested family of
compact convex subsets, with
$K_{m+1}$ contained in the interior of
$K_m$ for all $m \in \Z$.  Suppose that there
is some uniform constant $C$ such that the Hilbert
diameter of $K_{m+1}$ with respect to $K_m$ 
is less than $C$ for all $m$.  Then
$\bigcap K_m$ is a single point.
\end{corollary}

\startproof
We just use the upper bound
on diameter.  Lemma \ref{nest} implies
that the Euclidean diameter
of $K_{m+1}$ is at most $\lambda$ times the Euclidean
diameter of $K_m$, for some uniform $\lambda<1$.
\endproof

Now we turn to the triangulation $X$ associated to a
pentagram  spiral $P$.  We think of $X$ as
a union of solid triangles. 

Let $(A_0,B_0)$ be the seed generating $P$.
Define
\begin{equation}
(A_m,B_m)=T^{mn}(A_0,B_0).
\end{equation}
Let $K_m$ be the convex region bounded by the
polygon $A_m$.  By construction,
$K_{m+1}$ is contained in the interior of
$K_m$. 

By Theorem \ref{compact} and compactness, there is an
upper bound $C$ such that the Hilbert diameter of
$K_{m+1}$ with respect to $K_m$ is less than $C$,
independent of $m$. Here we are using the
projective naturality of the Hilbert metric.
 The point is that we are just
sampling a pre-compact subset of seeds in
${\cal C\/}(n,k)$.
Corollary \ref{shrink} now says that
\begin{equation}
\bigcap K_m
\end{equation}
is a single point.  But
$X$ contains $K_m-K_{m+1}$.  Hence
$X$ contains all points of $K_0$ except
for this one interesection point. 
This is what we wanted to prove.

Next we claim that $X$ is a triangulation of
$\R\P^2-L$, where $L$ is a single line.
This is equivalent to the statement that
the $\omega$ limit set of $P_-$ is a single
line $L$.  Here $P_-$ is the outward
portion of $P$.  

As is well known, projective dualities conjugate the
(suitably interpreted) pentagram map to its inverse.
Thus, the $\omega$-limit set of $P_-$ is dual to the
$\omega$-limit set of $P^*_+$, the inward part of
the dual pentagram spiral.  We have already proved that
the latter is a single point.  Hence, the former is
a single line.   

This completes the proof of Theorem \ref{omega}.

\subsection{Winding of Pentagram Spirals}

Let $P$ be a pentagram spiral, with distinguished
vertex $P_1$.  We translate the picture so that
$P \subset \R^2$ and so that the origin is the
limit point and $P_1$ is on the positive
$x$-axis, as shown in Figure 6.2.

\begin{center}
\resizebox{!}{4.5in}{\includegraphics{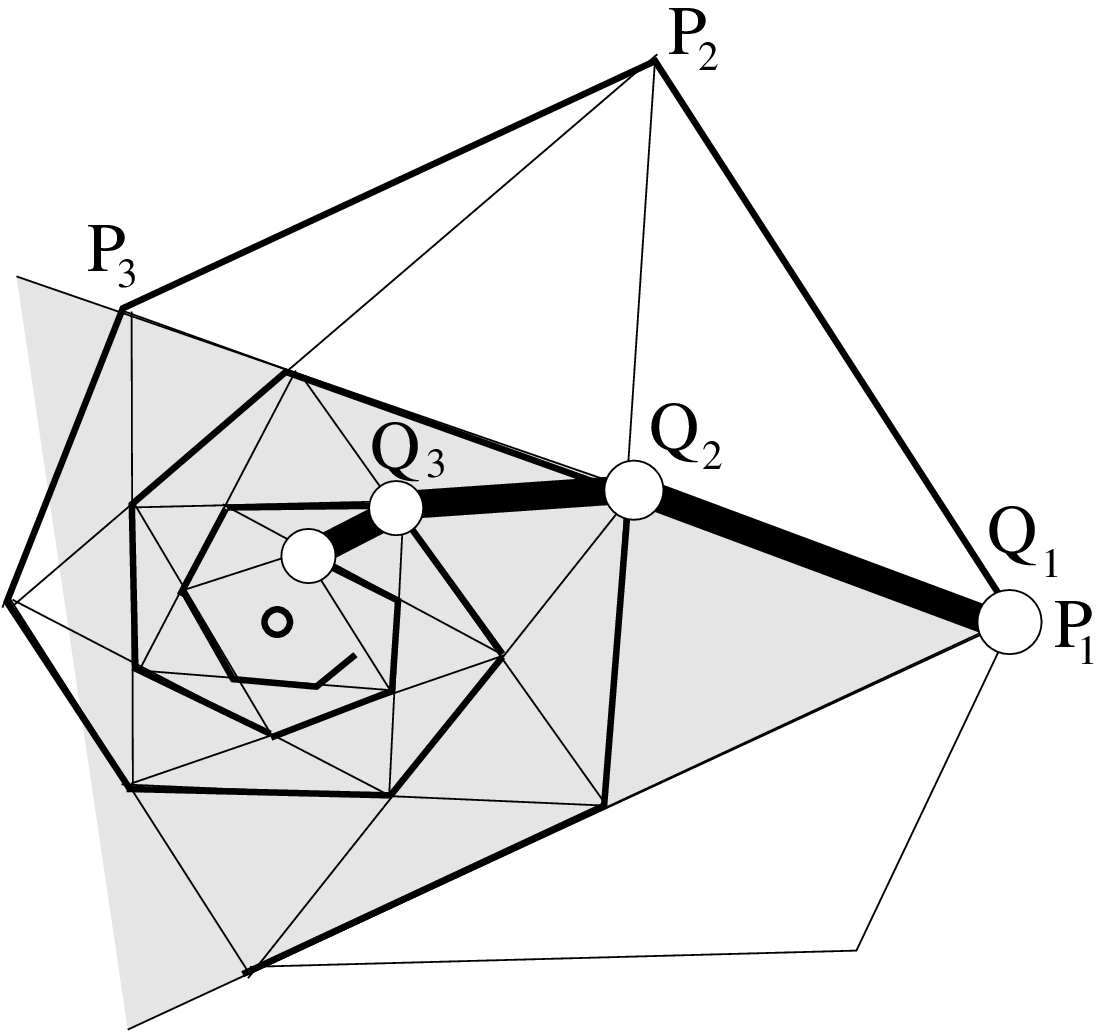}}
\newline
{\bf Figure 6.2:\/} The spine
\end{center}

The vertex $Q_1=P_1$ is the apex of a (shaded) cone
$C$ which contains the inward spiraling 
direction of $P$ starting at
$P_3$. Hence $C$
contains the limit point of $P$, namely the origin.
The top edge of $C$ leads to a point $Q_2$ lying
in the upper half plane.  The important point
here is that $\arg(Q_2)>\arg(Q_1)$.
We can repeat the same construction at $Q_2$. There
is a cone which contains the origin, and one of
the edges of this cone leads to a point $Q_3$ 
such that $\arg(Q_3)>\arg(Q_2)$.

Let $P_+$ be the inward spiraling direction of $P$.
Let $Q_+$ be the polygonal path connecting
the points $Q_1,Q_2,Q_3,...$ as shown in Figure 
6.2. 
 Both these paths converge to the origin in the 
forward direction.  
Since $Q_+$ contains infinitely many points and
each such point lies on one of the finitely many
spirals in the tiling, we see that some spiral contains
infinitely many point of $Q_+$.  But, in fact, by
inspecting the picture we see that whenever some
point of $Q_+$ intersects a spiral, the next point
of $Q_+$ intersects the next spiral.  Thus $P_+$ and
$Q_+$ intersect infintely often.

  Since the argument increases along
both $P_+$ and along $Q_+$, we see that
the argument increases by more than $2\pi$
along $P_+$ between every two intersections
with $Q_+$.  Hence, the argument along
$P_+$ increases without bound.
This proves Theorem \ref{wind}.
\newline
\newline
{\bf Remark:\/}
The proof of Theorem \ref{wind} is really quite
simple.  The only way it relies on the previous
material is that we would like to say that
the pentagram spiral $P$ really does have a single
limit point.

\newpage
\section{Experiments and Discussion}

\subsection{Periodicity}

With respect to specially
chosen labeling schemes, the pentagram map
is the identity on ${\cal C\/}(5)$ and
has period $2$ on ${\cal C\/}(6)$.
See [{\bf Sch1\/}].  Referring to the
pentagram spirals,
I discovered the following result computationally.

\begin{theorem}
\label{period}
The following is true.
\begin{itemize}
\item $T_{4,1}^2$ is the identity on
${\cal C\/}(4,1)$.
\item $T_{4,2}^2$ is the identity on
${\cal C\/}(4,2)$.
\item $T_{5,1}^8$ is the identity on
${\cal C\/}(5,1)$.
\end{itemize}
\end{theorem}

Theorem \ref{period} says that all the
pentagram spirals of this kind are
self-projective.
Theorem \ref{period} can be expressed
in terms of seeds and the seed map,
and so it only involves a finite number
of points and lines.  Thus, Theorem
\ref{period} can be established by a
finite calculation, similar in spirit
to what is done in [{\bf ST\/}].  I
have not yet made the rigorous calculations.
Some of the results in [{\bf ST\/}] have
conceptual proofs, and I wonder if there
are likewise conceptual proofs for the statements
in Theorem \ref{period}

It seems discussing how I discovered this.
My computer program allows the user to 
normalize the spirals so that the first
$4$ points are the vertices of the
unit square. One can then watch an animation
which shows the iteration of $T_{n,k}$.
I put in an option which allows the user
to choose a smallish integer $q$ and watch
the movie showing $T_{n,k}^q$.

For instance, for the parameter $(n,k)=(4,3)$ the
choice $q=18$ makes for a nice movie.  The point
is that, for a random choice of pentagram spiral
of type $(4,3)$, the $18$th power of the shift
map is fairly close to the identity, so an
animation of the map looks somewhat like a flow
to the naked eye. As another
example, for $(n,k)=(6,2)$ the choice $q=54$
often produces a nice movie.  I have found a
few of these values experimentally, but not many.
The reader can see all of this in action
on my program.

For the combinatorial type $(4,1)$, I noticed that
the image on the computer screen, for $q=1$ just
flickered back and forth.  When I put in $q=2$
the image was just stationary.  Likewise, for
the combinatorial type $(4,2)$, the choice
$q=2$ ``froze the movie'' and for the combinatorial
type $(5,1)$ the choice $q=8$ ``froze the movie''.
I would say that this is overwhelming experimental
support for Theorem \ref{period}.

\subsection{Asymptotic Shape}

We can combine Theorem \ref{period} with
some of our other results to get more
information about the special cases
$(4,1)$, $(4,2)$.

\begin{theorem}
Any PLC pentagram spiral of type $(4,1)$ of
$(4,2)$ is projectively equivalent to a
self-similar polygonal path.
\end{theorem}

\startproof
We can normalize by a projective
transformation so that our pentagram spirals
lie in the plane and have the origin as their
limit point.  If $P$ is a pentagram spiral as
in Theorem \ref{period}, then there is a
projective transformation $S$ so that
$S(P)=P$.  Necessarily $S$ preserves
the $\omega$-limit set of $P$. Hence
$S$ preserves the line at infinity and
fixes the origin.  That is, $S$ is a
linear transformation.

Now, the action of $S$ shifts the indices of
$P$ by $2$, or $8$ depending on the case.
In the cases $(4,1)$ and $(4,2)$ one can
argue that the path made from the
short diagonals of $P$ again winds infinitely
often around the origin.  But this means that
the orbits of $S$ wind infinitely often around the
origin.  But then $S$ is conjugate to a
similarity. So, in the cases $(4,1)$ and $(4,2)$,
there is a canonical normalization of the pentagram
spirals so that they are self-similar.
\endproof

The case $(4,1)$ yields a
$1$-parameter family of distinct shapes and
the case $(4,2)$ yields a $2$ parameter family
of shapes.   The argument above breaks down
in the case $(5,1)$, but experimentally the
result seems to be true.

 To describe some experimental
results along these lines, we need to introduce some
terminology. Let $\Gamma$ be an infinite polygonal
path which limits on the origin in one direction and
exits every compact set of the plane in the other
direction.  We call $\Gamma$ {\it quasi-logarithmic\/}
(or QL for short)
if there is some nontrivial homothety $D$ such that
the family of curves $\{D^n(\Gamma)\}$ is precompact
in the Hausdorff topology on shapes.  To say that
$\Gamma$ is QL is to say that $\Gamma$ is
only boundedly far from being invariant under a homothety.
For instance, a logarithmic spiral is quasi-logarithmic.

I mentioned above that my computer program normalizes
pentagram spirals so that $4$ distinguished vertices
are the vertices of a unit square.  I also programmed
things so that different normalizations are possible.
For instance, one can normalize by a homothety so that
the distinguished vertex is a point of the unit circle.
If a particular pentagram spiral is QL, then the movie
shown with this alternate normalization would ``go on
forever'' keeping more or less the same shape.

It seems that for most choices of $n$ and $k$, the
pentagram spirals of type $(n,k)$ are QL.  For
instance, when $n \leq 6$, only the case $(6,1)$
seems to produce pentagram spirals which are not
QL.   In general, it seems that the larger values of
$k$ produces pentagram spirals which are more likely
to be QL.  I would need to do more experiments before making
a definitive conjecture on this, so let me just
pose this as a question:
\newline
\newline
{\bf Question:\/} Are there values $(n,k)$ such that
every PLC pentagram spiral is QL?  If so, which values?

\subsection{Logarithmic Pentagram Spirals}

For each choice of $(n,k)$ there is a self-similar
pentagram spiral $P$ which I call the
{\it logarithmic pentagram spiral\/} (or
LPS for short.)  The vertices of $P$ lie
in a logarithmic spiral and the edges of $P$
are inscribed in a rotated copy of the  logarithmic
spiral.
The pentagram map carres $P$ to a rotated copy
of $P$.  

The LPS can be
normalized so that it has
vertices $\{z^n|\ n \in \Z\}$, where
\begin{equation}
\label{spiralX}
|z|<1, \hskip 20 pt
\frac{2\pi}{n+k}<\arg(z)<
\frac{2\pi}{n}, \hskip 20 pt
(z+\overline z)^k = z^{n+k}(1+\overline z)^k.
\end{equation}
These equations come from the observation that
the intersection of the line through $1$ and $z^2$
with the line through $z$ and $z^3$ is
$$w=\frac{z(z+\overline z)}{1+\overline z},$$
and that the combinatorics of the spiral dictate
that $w^k=z^{n+2k}.$
I had a lot of trouble solving Equation \ref{spiralX} on
Mathematica, but I will describe a different way to
draw very close approximations to the LPS.

Say that the seed $P=(A,B)$ is
{\it normalized\/} if
\begin{equation}
A_1=(0,0), \hskip 20 pt
A_2=(1,0), \hskip 20 pt
A_3=(1,1), \hskip 20 pt
A_4=(0,1).
\end{equation}
That is, the first $4$ vertices of $P$ are the vertices
of the unit square.  Each point of ${\cal C\/}(n,k)$ has
a unique normalized representative.
A normalized representative consists of a convex $n$-gon
and $k$ additional points.  The locations of these
points can be described by the quantities
\begin{equation}
d_i=\frac{\|A_{i}-B_i\|}{\|A_i-A_{i+1}\|} \in (0,1), \hskip 30 pt
i=(n-k+1),...,n.
\end{equation}
The point $A_{n+1}$ is interpreted as $A_1$.

Given some point $P=(A,B) \in {\cal C\/}(n,k)$ and some
integer $j$, we define $P^j=(A^j,B^j)$ as the normalized
version of $T_{n,k}^j(P)$.  Given some integer $m$,
define $\Theta_m(P)=(A',B')$, where $A'$ is the pointwise
average of the polygons $A^0,....,A^{m-1}$ and
the number $d_j'$ is the average of the
corresponding numbers for $P^0,...,P^{m-1}$.
In practice, $\Theta_m$ seems to act as a
kind of contraction
mapping on ${\cal C\/}(n,k)$, and the fixed point is
the logarithmic pentagram spiral. 

So, one can choose
some smallish $m$ and iterate $\Theta_m$ several
times.  This produces a point in 
${\cal C\/}(n,k)$ very close to the point
representing the logarithmic pentagram spiral.
Given a point in ${\cal C\/}(n,k)$ representing
the logarithmic pentagram spiral, one can then
apply projective transformations to find
better normalizations.  The cheapest way to
do this is just to apply $T_{n,k}$ many times
and then dilate the picture.
This method produced the picture shown
in Figure 7.1.

\begin{center}
\resizebox{!}{3.5in}{\includegraphics{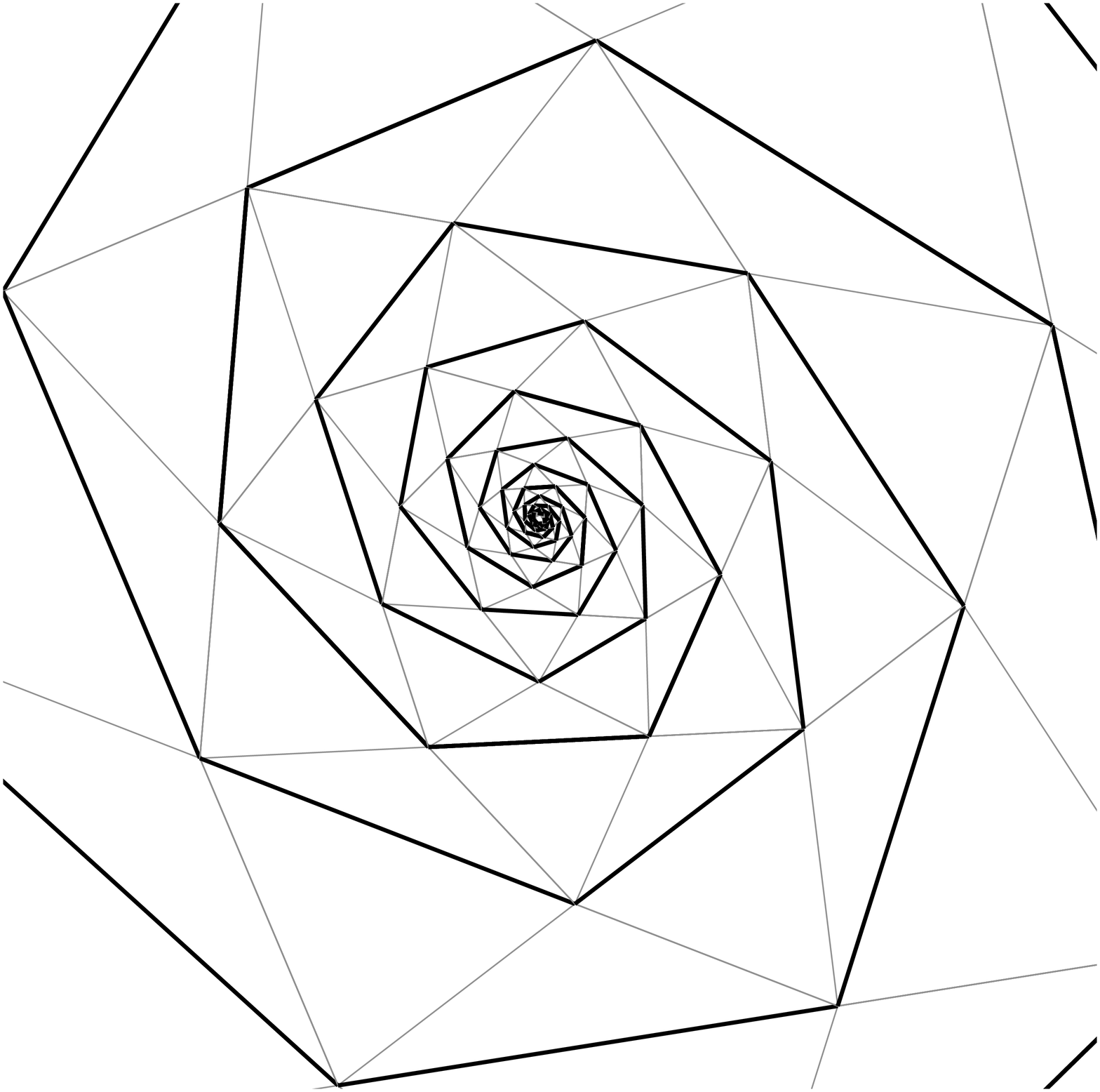}}
\newline
{\bf Figure 7.1:\/} Approximation to the
LPS of type $(5,3)$.
\end{center}

The logarithmic pentagram spiral is a natural origin
for the space ${\cal C\/}(n,k)$ just as the
projective class of the regular $n$-gon is a
natural origin for the space ${\cal C\/}(n)$.
Computer experiments suggest the following
conjecture

\begin{conjecture}
For any $(n,k)$, the invariant
$Z$ is uniquely maximized, and has a unique critical
point, at the point representating the logarithmic
pentagram spiral.
\end{conjecture}

There is an analogous conjecture for the space
${\cal C\/}(n)$.  Corollary 1.2 in my paper
[{\bf Sch4\/}] proves that (the analogue of)
$Z$ is maximized at the regular polygons when
restricted to the subspace of polygons
which are inscribed in circles.

\subsection{Twisted Pentagram Spirals}

The notion of a twisted polygon has been very useful
in the study of the ordinary pentagram map.
See [{\bf Sch3\/}], [{\bf OST1\/}], and
[{\bf OST2\/}].
A {\it twisted $n$-gon\/} is a map
$\phi: \Z \to \R\P^2$ which intertwines
translation by $n$ with a projective transformation
$M$.  That is,
\begin{equation}
\label{mono}
\phi \circ \mu=M \circ \phi(k), \hskip 30 pt
\end{equation}
Here $\mu$ is translation by $n$.  That is,
$\mu(k)=k+n$.
The transformation $M$ is the {\it monodromy\/}
of the twisted polygon.

The pentagram map acts on twisted $n$-gons
and again commutes with projective transformations.
One can define $2n$ flag invariants of a twisted
$n$-gon just as for an ordinary ones.  The
flag invariants of twisted $n$-gon and its
iterates under the pentagram map naturally
give rise to a labeling of $\T$, as in
\S \ref{tiling2}. This labeling is an element
of ${\cal T\/}(n,0)$.

Now we want to do the same kind of thing for
pentagram spirals.
First of all, we fix a {\it model\/} tiling
generated by a particular element of ${\cal C\/}(n,k)$,
say the logarithmic pentagram spiral.
(The model is just used for combinatorial purposes.)
Next, we create a locally
identical tiling of the universal cover $\widetilde X$ of
$\R^2-(0,0)$.  That is, we simply pull back
the tiling to $\widetilde X$.
Let $\widetilde T$ denote the tiling of
the universal cover.

The space $\widetilde X$ is 
homeomorphic to the plane, but it
has an exotic projective structure on it.
A straight line in this universal cover
is something which projects to a straight line.
There is a $\Z$-action on $\widetilde X$,
namely the deck group.  The deck group
acts so as to carry lines to lines.
Let $\mu$ be a generator of the deck group.
We have $\mu(\widetilde T)=\widetilde T$.

A {\it twisted pentagram spiral\/} is an
adapted map $\phi$ satisfying
Equation \ref{mono} with respect to the
deck group generator $\mu$.
By {\it adapted map\/} I mean
that $\phi$ is a homeomorphism when restricted
to each (solid) triangle of $\widetilde T$ and that
$\phi$ carries each line segment in
the $1$-skeleton of $\widetilde T$ to a
straight line segment in $\R\P^2$.
Note that some of these line segments consist
of $3$ consecutive edges of triangles.

An ordinary pentagram spiral is simply a
twisted pentagram spiral having monodromy
the identity.  Moreover, the above
definition reduces to the
usual definition in the case of
closed polygons.  It is merely the
original definition rephrased in terms
of the triangulation produced by the
pentagram map.

Some readers might not like the above geometric
definition of a twisted pentagram spiral, so let
me describe things algebraically.
A projective equivalence class of
A twisted pentagram spiral is nothing
more than an element
of ${\cal T\/}(n,k)$.  Starting
with an element of ${\cal T\/}(n,k)$, one
can start building a network of line segments
in the projective plane, such that the
corresponding flag invariants give the labels
of ${\cal T\/}(n,k)$.  As one develops
the picture going ``all the way around'' the
cylinder $\R^2/V_{n,k}$, one might observe
that the configuration in the projective plane
does not close up.  The
failure of the picture to close up is
encoded in the monodromy.

\subsection{Integrability}

Computer experiments suggest the following
conjecture.

\begin{conjecture}
\label{integrability}
The map
$T_{n,k}$ acting on ${\cal C\/}(n,k)$ is
a discrete totally integrable system.
\end{conjecture}

What I mean is that ${\cal C\/}(n,k)$
should have a singular foliation by tori, each
equipped with a flat structure, such that
each orbit of $T_{n,k}$ is contained in
a finite union of tori.  Moreover, the
restriction of a suitable power of $T_{n,k}$
to each torus is a translation 
relative to the canonical flat structure.
Such a structure would arise if
${\cal C\/}(n,k)$ had an invariant
Poisson structure and sufficiently
many commuting invariant functions.
This how the torus foliation arises in
[{\bf OST1\/}] and [{\bf OST2\/}] for
the pentagram map.

I will describe, to some extent, the
computer experiments which lead to this conjecture.
The interested reader can do the
experiments themselves using my program.

Let $T$ be a PLC pentagram spiral, and
let $c(T)$ denote the limit point of the
inward spiraling direction of $T$.
If we normalize $T$ so that the first $4$
vertices are the vertices of the unit square
$Q$, then the point $c(T) \in \R^2$ is a canonical
point associated to $T$.  Assuming that $T$ is
a pentagram spiral of type $(n,k)$, we define
\begin{equation}
c_m=c(T_{n,k}^m(T)).
\end{equation}
That is, $c_m$ is the limit point of the $m$th iterate
of $T$ under the shift map.  
My computer program allows the reader to view
the sequence $\{c_m\}$.  

For instance, the points
$\{c_{m}\}$ seem to lie on a union of two (generically) smooth
curve when $(n,k)=(5,2)$.  When we consider the thinner sequence
$\{c_{18n}\}$ in this case, we see points appearing in order
on a smooth curve.  In other words, the movie we produce
for the choices $(n,k)=(5,2)$ and $q=18$ shows the limit
point gently rotating around a smooth curve.  Again,
we encourage the reader to download the program, so that
he or she can see this in action.

The space ${\cal C\/}(5,2)$ is $4$-dimensional. The
experiments above suggest that ${\cal C\/}(5,2)$ is
foliated by invariant loops and so that
$R_{5,2}^2$ preserves each loop in the foliation
and acts there as a (typically) irrational rotation.
As is the case with the pentagram map, one would
describe this situation as ``mildly hyperintegrable'':
The completely integrable situation would predict
invariant $2$-tori.

As $n$ and $k$ increase, it is harder to see that the
sequence $\{c_m\}$ is the projection of a sequence of
curves lying on a finite union of tori. However, for
smallish values of $m$ and $k$, one still gets a sense
that this is the case.

\subsection{Monodromy Invariants}

The first step to proving the integrability
conjecture is to find the integrals (i.e. invariants)
of the map $T_{n,k}$.  
For the pentagram map, these invariants by
now have many constructions. I will describe
the original way I thought of them.
Let $M$ be the monodromy of a twisted $N$-gon.
If we replace the twisted $N$-gon by a
projectively equivalent one, the
monodromy $M$ is replaced by a conjugate.
However, the two quantities
\begin{equation}
\label{trace}
\frac{{\rm Tr\/}(M)}{\det^{1/3}(M)}, \hskip 50 pt
\frac{{\rm Tr\/}(M^*)}{\det^{1/3}(M^*)}.
\end{equation}
only depend on the projective equivalence class.
In Equation \ref{trace} we think of $M$ and
$M^*$ as matrices representing the action
of the monodromy on the projective plane
and on the dual projective plane respectively.

The quantities in Equation \ref{trace}
 are rational functions of the
flag invariants.  There is a certain natural
weighting of the monomials in these rational
functions, and the homogeneous parts with
respect to this weighting are the monodromy
invariants. See [{\bf Sch3\/}], [{\bf OST1\/}]
and [{\bf OST2\/}] for details about this.

The special weighting can be described as follows.
If we take any element of ${\cal T\/}(n,k)$ we
can multiply all the forward slanting edges by $s$
and all the backward slanting edges by $1/s$.  This
produces a new element of ${\cal T\/}(n,k)$. All
the compatibility equations hold for the new
equations.
Every paper which has discussed
the integrability of the pentagram map (and
its generalizations) uses this scaling in a
crucial way.

Now, the monodromy of a pentagram spiral can
be computed using a path which only encounters
finitely many edges of the tiling $\widetilde T$.
This, it would seem that the quantities in
Equation \ref{trace} would also be rational
functions in finitely many of the flag
invariants. The scaling mentioned above
works just fine here.  So, the weighted homogeneous
parts ought to be invariants of the map
$T_{n,k}$ on the larger space 
${\cal T\/}(n,k)$, which we might as well interpret
as the space of twisted pentagram spirals of
type $(n,k)$.  Restricting these invariants
to the subset ${\cal C\/}(n,k)$ we would
get invariants for the shift map $T_{n,k}$.

I tried to compute the monodromy invariants for
the very modest case $(n,k)=(4,1)$ and I arrived
at depressingly complicated expressions.  This
makes me somewhat pessimistic that one could
arrive at crisp formulas like
Equation \ref{trace} in the general case.
It seems to me that the calculation will have
to wait either for a more determined 
experimenter or for a better coordinate system.

My guess is that the Poisson bracket of
[{\bf OST1\/}] and [{\bf OST2\/}] will
generalize to the case of pentagram spirals
as well.  But I will leave this discussion
for a later paper.

\newpage

\section{References}

[{\bf FM}]
V. Fock, A. Marshakov,
{\it Integrable systems, clusters,  dimers and loop groups},
preprint, 2013
\newline \newline
[{\bf GK\/}] A. B. Goncharov and R. Kenyon,
{\it Dimers and Cluster Integrable Systems\/}, 
preprint, arXiv 1107.5588, 2011.
\newline \newline
[{\bf GSTV}]
M. Gekhtman, M. Shapiro, S.Tabachnikov, A. Vainshtein,
{\it Higher pentagram maps, weighted directed networks, and cluster dynamics}, 
Electron. Res. Announc. Math. Sci. {\bf 19} 21012, 1--17
\newline \newline
[{\bf Gli}]
M. Glick, 
{\it The pentagram map and $Y$-patterns}, 
Adv. Math. {\bf 227}, 2012,  1019--1045.
\newline \newline
[{\bf Gli1}]
M. Glick, 
{\it The pentagram map and $Y$-patterns},
23rd Int. Conf. on Formal Power Series and Alg. Combinatorics (FPSAC 2011), 
399--410.
\newline\newline
[{\bf KDif\/}], R. Kedem and P. DiFrancesco, {\it $T$-Systems with boundaries from network solutions\/}, preprint, arXiv 1208.4333, 2012
\newline \newline
[{\bf KS}]
B. Khesin, F. Soloviev
{\it Integrability of higher pentagram maps},
Mathem. Annalen. (to appear) 2013
\newline \newline
[{\bf MB1}]
G. Mari Beffa,
{\it On Generalizations of the Pentagram Map: Discretizations of AGD Flows\/},
arXiv:1303.5047, 2013
\newline \newline
[{\bf MB2}]
G. Mari Beffa,
{\it On integrable generalizations of the pentagram map} \newline
arXiv:1303.4295, 2013
\newline \newline
[{\bf Mot\/}]
Th. Motzkin, 
{\it The pentagon in the projective plane, with a comment on NapierÕs rule},
Bull. Amer. Math. Soc. {\bf 52}, 1945, 985--989.
\newline \newline
[{\bf OST}]
V. Ovsienko, R. Schwartz, S. Tabachnikov, {\it Quasiperiodic motion for the pentagram map}, 
 Electron. Res. Announc. Math. Sci. {\bf 16} ,2009, 1--8.
\newline \newline
[{\bf OST1}]
V. Ovsienko, R. Schwartz, S. Tabachnikov, 
{\it The pentagram map: A discrete integrable system}, 
Comm. Math. Phys. {\bf 299}, 2010, 409--446.
\newline \newline
[{\bf OST2}]
V. Ovsienko, R. Schwartz, S. Tabachnikov, 
{\it Liouville-Arnold integrability of the pentagram map on closed polygons},
to appear in Duke Math. J.
\newline \newline
[{\bf Sch1}]
R. Schwartz,
{\it The pentagram map},
Experiment. Math. {\bf 1}, 1992, 71--81.
\newline \newline
[{\bf Sch2}]
R. Schwartz,
{\it The pentagram map is recurrent},
Experiment. Math. {\bf 10}, 2001, 519--528.
\newline \newline
[{\bf Sch3}]
R. Schwartz,
{\it Discrete monodromy, pentagrams, and the method of condensation},
J. of Fixed Point Theory and Appl. {\bf 3}, 2008, 379--409.
\newline \newline
[{\bf Sch4}]
R. Schwartz,
{\it A Conformal Averaging Process on the Circle\/}
Geom. Dedicata., 117.1, 2006.
\newline \newline
[{\bf Sol}]
F. Soloviev
{\it Integrability of the Pentagram Map}, to appear in Duke Math J.
\newline \newline
[{\bf ST\/}]
R. Schwartz, S. Tabachnikov,
{\it Elementary surprises in projective geometry},
Math. Intelligencer {\bf 32}, 2010, 31--34.

\end{document}